\documentclass[12pt]{article}

\usepackage{amsmath,amsfonts,amssymb,amsthm,latexsym,amscd}
\newcommand{\var}{\mbox{var}}
        
\newlength{\figboxwidth}             
\newcommand{\infinity}{\infty}
\newcommand{\integral}{\int}
\newcommand{\F}{{\mathcal F}}

\newcommand{\cross}{\times}

\newcommand{\ignore}[1]{}  

\newtheorem{theorem}{Theorem}[section]

\newtheorem{lemma}[theorem]{Lemma}

\newtheorem{corollary}[theorem]{Corollary}
\newtheorem{definition}[theorem]{Definition}
\newtheorem{question}[theorem]{Question}
\newtheorem{remark}[theorem]{Remark}
\newtheorem{claim}[theorem]{Claim}

\mathchardef\GG="321D
\usepackage{color}

\definecolor{Red}{rgb}{1,0,0}
\definecolor{Blue}{rgb}{0,0,1}
\definecolor{Olive}{rgb}{0.41,0.55,0.13}
\definecolor{Yarok}{rgb}{0,0.5,0}
\definecolor{Green}{rgb}{0,1,0}
\definecolor{MGreen}{rgb}{0,0.8,0}
\definecolor{DGreen}{rgb}{0,0.55,0}
\definecolor{Yellow}{rgb}{1,1,0}
\definecolor{Cyan}{rgb}{0,1,1}
\definecolor{Magenta}{rgb}{1,0,1}
\definecolor{Orange}{rgb}{1,.5,0}
\definecolor{Violet}{rgb}{.5,0,.5}
\definecolor{Purple}{rgb}{.75,0,.25}
\definecolor{Brown}{rgb}{.75,.5,.25}
\definecolor{Grey}{rgb}{.5,.5,.5}

\newcommand{\suml}{\sum\limits}

\newcommand{\R}{\mathbb{R}}
\newcommand{\Z}{\mathbb{Z}}
\newcommand{\halb}{\frac{1}{2}}
\newcommand{\prob}{{\bf P}}
\newcommand{\E}{{\bf E}}
\newcommand{\C}{{\cal C}}
\newcommand{\en}{{\cal E}}
\newcommand{\idp}{{\cal I}}
\newcommand{\ef}{{\rm eff}}
\newcommand{\borel}{{\mathcal B}}
\newcommand{\fernandez}{$\text{Fern}\acute{\text{a}}\text{ndez }$}
\newcommand{\re}{\text{  }}
\theoremstyle{definition}
\theoremstyle{lemma}

\title{Transience, Recurrence and Critical Behavior for Long-Range Percolation}
\author{Noam Berger\footnote{Research  partially supported by NSF grant
 \#DMS-9803597 and
by a US-Israel BSF grant.
}\thanks{Part of the research was done while the author was at the Hebrew University of Jerusalem.}
\\The University of California at Berkeley}
\date{August 2001}

\begin{document}
\maketitle

\begin{abstract}
We study the behavior of the random walk on the infinite cluster
of independent long-range percolation in dimensions $d=1,2$, where $x$ and $y$ are
connected with probability $\sim\beta/\|x-y\|^{-s}$.
We show that if $d<s<2d$ then the
walk is transient, and if $s\geq 2d$, then the walk is recurrent. The proof of
transience is based on a renormalization argument. As a corollary of this
renormalization argument, we get that for every dimension $d\geq 1$,
if $d<s<2d$, then
there is no infinite cluster at criticality. This result is extended to the
free random cluster model.
A second corollary is that when $d\geq 2$ and $d<s<2d$ we can erase
all long enough bonds and still have an infinite cluster.
The proof of recurrence in two dimensions is based on general stability
results for recurrence in random electrical networks. In particular, we show
that i.i.d. conductances on a recurrent graph of bounded degree yield a
recurrent electrical network.
\end{abstract}

\section{Introduction}
\subsection{background}\label{background}
Long-range percolation (introduced by Schulman in 1983 \cite{schul}) is a percolation
model on the integer lattice $\Z^d$ in which every
two vertices can be connected by a bond.
The probability of the bond between two vertices to be open depends on
the distance between the vertices.
The models that were studied the
most are models in which the probability of a bond to be open decays polynomially with
its length.

\subsection{The model - definitions and known results}
Let $\{P_k\}_{k\in\Z^d}$ be s.t. $0\leq P_k=P_{-k} < 1$ for every $k\in\Z^d$.
We consider the following percolation model on $\Z^d$: for every $u$ and $v$ in $\Z^d$,
the bond connecting $u$ and $v$ is open with probability $P_{u-v}$. The different bonds
are independent of each other.

\begin{definition}
For a function $f:\Z^d\to\R$, we say that $\{P_k\}$ is {\em asymptotic to} $f$ if
\begin{equation*}
\lim_{\|k\|\to\infty}\frac{P_k}{f}=1.
\end{equation*}
We denote it by $P_k\sim f(k)$
\end{definition}

Since the model is shift invariant and ergodic, the event that an infinite cluster exists
is a zero--one event. We say that $\{P_k\}$ is {\em percolating} if a.s. there exists an
infinite cluster.

We consider systems for which $P_k\sim \beta \|k\|_1^{-s}$ for certain $s$ and $\beta$.
The following facts are trivial.
\begin{itemize}
\item
If $s\leq d$, then $\sum_k{P_k}=\infty$. Therefore, By the Borel Cantelli Lemma,
every vertex is connected to infinitely many other vertices. Thus, there exists an infinite
cluster.
\item
If $\sum_k{P_k}\leq 1$ then by domination by a (sub)-critical Galton-Watson tree
there is no infinite cluster. Therefore,
for every $s>d$ and $\beta$ one can find a set $\{P_k\}$ s.t. $P_k\sim \beta {\|k\|_1^{-s}}$
and s.t. there is no infinite cluster.
\end{itemize}

In \cite{schul}, Schulman proved that if $d=1$ and $s>2$, then there is no infinite cluster.
Newman and Schulman (\cite{NS}) and Aizenman and Newman (\cite{AN}) proved,
among other results, the following:
\begin{theorem}
(A)
If $d=1$, { } $1<s<2$, and $P_k\sim \beta |k|^{-s}$ for some $\beta>0$, then there exists a
$\{P'_k\}$ s.t. $P'_k=P_k$ for every $k\geq 2$, { } $P'_1<1$ and $\{P'_k\}$ is percolating.
I.e., if $1<s<2$ then by increasing $P_1$ one can make the system percolating.
\medskip
\\(B)
If $d=1$, { } $s=2$, { } $\beta>1$, and $P_k\sim \beta |k|^{-s}$, then there exists a
$\{P'_k\}$ s.t. $P'_k=P_k$ for every $k\geq 2$, { } $P'_1<1$ and $\{P'_k\}$ is
percolating.
\medskip
\\(C)
If $d=1$, { } $s=2$, { } $\beta\leq 1$, and $P_k\sim \beta |k|^{-s}$ then $\{P_k\}$ is not
percolating.
\end{theorem}
These results show the existence of a phase transition for $d=1$, { } $1<s<2$ and
$\beta>0$, and for $d=1$, { } $s=2$ and $\beta>1$.

When considering $\Z^d$ for $d>1$, the picture is simpler. The following fact is a trivial
implication of the existence of infinite clusters for nearest-neighbor percolation:
\begin{itemize}
\item
If $d>1$, { } $s>d$ and $P_k\sim \beta {\|k\|_1^{-s}}$ for some $\beta>0$, then there exists a
percolating
$\{P'_k\}$ s.t. $P'_k=P_k$ for every $\|k\|_1\geq 2$ and $P'_k<1$ for every $k$ whose norm is
$1$.
\end{itemize}
If $d>1$, then for any $s>d$ and $\beta>0$ we may obtain a transition between the
phases of existence
and non-existence of an infinite cluster by only changing $\{P_k|k\in A\}$ for a finite
set $A$.

In \cite{uniq}, Gandolfi, Keane and Newman proved a general uniqueness theorem.
A special case of it is the following theorem:
\begin{theorem}\label{GKN}
If $\{P_k\}_{k\in\Z^d}$ is percolating and for every $k\in\Z^d$ there exist $n$ and
$k_1,...,k_n$ s.t. $k=k_1+k_2+...+k_n$ and $P_{k_i}>0$ for all $1\leq i \leq n$ then
a.s. the infinite cluster is unique.

In particular,
If $\{P_k\}_{k\in\Z^d}$ is percolating and $P_k\sim \beta {\|k\|_1^{-s}}$ for some $s$ and
$\beta>0$, then a.s. the infinite cluster is unique.
\end{theorem}

\subsection{Goals}
Random walks on percolation clusters have been studied intensively in
recent years. In \cite{GKZ}, Grimmett, Kesten and Zhang showed that a
supercritical percolation in $\Z^d$ is transient for all $d\geq 3$.
See also \cite{BPP}, \cite{elch} and \cite{us}.

The problem discussed in this paper,
suggested by Itai Benjamini, was to determine when a random walk on the long-range
percolation cluster is transient. In \cite{Jesper}, Jespersen and Blumen worked on
a model which is quite similar to the long-range percolation on $\Z$, and they
predict that when $s<2$ the random walk is transient, and when $s=2$ it is recurrent.

\subsection{Behavior of the random walk}\label{rwalk}
The main theorem proved here is:
\begin{theorem}\label{main}
(I) Consider long-range percolation on $\Z$  with parameters
$P_k\sim\beta |k|^{-s}$ such that a.s. there is an infinite cluster.
If $1<s<2$ then
the infinite cluster is
transient. If $s=2$, then the infinite cluster is recurrent.
\\(II) Let $\{P_k\}_{k\in\Z^2}$ be percolating for $\Z^2$
such that $P_{k}\sim\beta \|k\|_1^{-s}$.
If $2<s<4$ then the infinite cluster
is transient. If $s\geq 4$, then the infinite cluster is recurrent.
\end{theorem}

In Section \ref{trsprf}, we prove
the transience for the one-dimensional case where $1<s<2$ and
for the two-dimensional case where $2<s<4$.
Actually, we prove more - we show that for every $q>1$ there is a
flow on the infinite cluster with finite $q$-energy, where the $q$-energy of a
flow $f$ is defined as
\begin{equation}\label{energy}
\en_q(f)=\sum_e{f(e)^q}.
\end{equation}
It is well known that finite $2$-energy is equivalent to
transience of the random walk (see e.g. \cite{yuval}, section 9),
so the existence of such flows is indeed a
generalization of the transience result (See also \cite{LP}, \cite{HoMo} and \cite{us}).

In Section \ref{recprf} we prove the recurrence for the one-dimensional case
with $s=2$ and for the two-dimensional case with $s\geq 4$.

\subsection{Critical behavior}\label{critbehav}
As a corollary of the main renormalization lemma, we prove the following
theorem, which applies to every dimension:
\begin{theorem}\label{criti}
Let $d\geq1$ and let $\{P_k\}_{k\in\Z^d}$ be probabilities such that
$P_k\sim\beta \|k\|_1^{-s}$. Assume that $d<s<2d$.
Then, if $\{P_k\}$ is
percolating then it is not critical, i.e. there exists an $\epsilon>0$ such
that the sequence $\{P'_k=(1-\epsilon)P_k\}$ is also percolating.
\end{theorem}

In \cite{hasl}, Hara and Slade proved, among other results, that for dimension $d\geq6$ and an
exponential decay of the probabilities, there is no infinite cluster at criticality.

It is of interest to compare Theorem \ref{criti}
with the results of Aizenman and Newman (\cite{AN}), that show that for $d=1$
and $s=2$, a.s. there exists an infinite cluster at criticality. In \cite{ACCN},
Aizenman, Chayes, Chayes and Newman showed the same result for the Ising model -
They showed that if $s=2$, then at the critical temperature there is a
non-zero magnetization.

The technique that is used to prove Theorem \ref{criti} is used in Section
\ref{frcm} to prove the analogous result for
the infinite volume limit of the free random cluster model, and to get:
\begin{theorem}\label{fr_int}
Let $\{P_k\}$ be a sequence of nonnegative numbers such that
$P_k\sim \|k\|_1^{-s}$ ($d<s<2d$) and let $\beta>0$. Consider the infinite volume limit of the
free random cluster model with probabilities $1-e^{-\beta P_k}$ and with
$q\geq 1$ states.
Then, at the critical inverse temperature
\begin{equation*}
\beta_c=\inf(\beta |\text{ a.s. there exists an infinite cluster})
\end{equation*}
there is no infinite cluster.
\end{theorem}

However, this
technique fails to prove this result for the wired measure, so in the wired
case the question is still open. A partial answer for the case $s\leq\frac{3}{2}d$ is given
by Aizenman and \fernandez in \cite{AiFe}. Consider the
Ising model with $s\leq\frac{3}{2}d$ when
the interactions obey the {\em reflection positivity}
condition (which is defined there). Denote by $M(\beta)$ the magnetization at
inverse temperature $\beta$. Consider the critical exponent $\hat{\beta}$ such that
\begin{equation*}
M(\beta)\sim |\beta-\beta_c|^{\hat{\beta}}
\end{equation*}
for $\beta$ near the critical value $\beta_c$. They proved that (under the above assumptions)
$\hat{\beta}$ (as well as other critical exponents) exists and they showed that
$\hat{\beta} = \halb$.

A corollary of Theorem \ref{fr_int} is
\begin{corollary}\label{ising_extr}
Let $\{P_k\}_{k\in\Z^d}$ be nonnegative numbers s.t. $P_k=P_{-k}$ for every $k$
and s.t. $P_k\sim \|k\|_1^{-s}$ ($d<s<2d$).
Consider the Potts model with $q$ states on $\Z^d$, s.t. the interaction between
$v$ and $u$ is  $P_{v-u}$. At the critical temperature,
the free measure is extremal.
\end{corollary}

Another consequence of the renormalization lemma is the following:
\begin{theorem}\label{jeff}
Let $d>1$ and let $\{P_k\}_{k\in \Z^d}$ be probabilities s.t.
$P_k\sim \beta \|k\|_1^{-s}$ for
some $s<2d$. Assume that the independent percolation model with $\{P_k\}$ has
a.s. an infinite cluster. Then there exists $N$ s.t. the independent
percolation model with probabilities
\begin{equation*}
P'_k = \left\{
        \begin{array}{ll}
                P_k             &\|k\|_1<N\\
                0               &\|k\|_1\geq N
        \end{array}
        \right.
\end{equation*}
also has, a.s., an infinite cluster.
\end{theorem}
In \cite{steif}, Meester and Steif prove the analogous result for supercritical
arrays of exponentially
decaying probabilities. It is still unknown whether the same statement is true
for probabilities that decay faster than $\|k\|_1^{-s}$ ($s<2d$) and slower than
exponentially.


\subsection{Random electrical networks}
The proof of
recurrence for the two-dimensional case involves some calculations on random
electrical networks. In Section \ref{elnet} we study such networks,
and prove
stability results for their recurrence. One of our goals in that section is:
\begin{theorem}\label{iid}
Let $G$ be a recurrent graph with bounded degree. Assign i.i.d
conductances on the edges of $G$. Then, a.s., the resulting electrical
network is recurrent.
\end{theorem}

In \cite{triid} Pemantle and Peres studied the analogous question for the transient case,
i.e.
under what conditions i.i.d. conductances on a transient graph would preserve
its transience. They proved that it occurs if and only if there exists $p<1$ s.t.
an infinite cluster for (nearest-neighbor) percolation with parameter $p$ is transient.

Comparing the results indicates that recurrence is more stable than transience for
this type of perturbation.

Section \ref{elnet} is self-contained, i.e. it does not use any of the results proved
in other sections.

\section{The transience proof}\label{trsprf}

In this section we give the proof that the $d$-dimensional long-range
percolation cluster, with $d<s<2d$, is a transient graph.
Our methods use the idea of iterated renormalization for
long-range percolation that was introduced in \cite{NS}, where it was used in order
to prove the following theorem:

\begin{theorem}[Newman and Schulman, 1986]\label{thm:ns}
Let $1<s<2$ be fixed, and consider  an independent  one-dimensional
percolation model such that the bond between $i$ and $j$ is open with
probability $P_{i-j}=\eta_s(\beta,|i-j|)$, where
\begin{equation}\label{defeta}
\eta_s(\beta,k)=1-\exp(-{\beta |k|^{-s}}),
\end{equation}
and each vertex is alive with probability $\lambda\leq 1$. Then for $\lambda$
sufficiently close to $1$
and $\beta$ large enough, there exists, a.s., an infinite cluster.
\end{theorem}
%
%

In order to prove our results, we need the following definition and
the following two renormalization lemmas:

\begin{definition}
We say
that the cubes
$\C_1=v_1+[0,N-1]^d$ and $C_2=v_2+[0,N-1]^d$ are {\em $k$ cubes
away from each other} if $\|v_1-v_2\|_1=Nk$.
\end{definition}
We will always use the notion of two cubes being $k$ cubes away from each other for
pairs of cubes of the same size that are aligned on the same grid.

\begin{lemma}\label{normal1}
Let $\{P_k\}_{k\in\Z^d}$ be such that $P_k=P_{-k}$ for every $k$, such that
$P_k>0$ for every $k\in\Z^d\setminus\{0\}$,
and such that
there exists $d<s<2d$ s.t.
\begin{equation}\label{posliminf}
\liminf_{\|k\|_1\to\infty}\frac{P_k}{\|k\|_1^{-s}}>0.
\end{equation}
Assume that the percolation model on $\Z^d$ with probabilities $\{P_k\}$ has,
a.s., an infinite cluster.

Then, for every $\epsilon>0$ and $\rho$ there exists $N$ such that with
probability bigger than $1-\epsilon$, inside the cube $[0,N-1]^d$
there exists an open cluster that contains at least $\rho N^{\frac{s}{2}}$ vertices.
\end{lemma}

Lemma \ref{normal1} shows that most of the cubes contain big clusters. We also
want to estimate the probability that the clusters in two different cubes are connected
to each other.

\begin{lemma}\label{connectprob}
Let $\{P_k\}_{k\in\Z^d}$ be such that $P_k=P_{-k}$ for every $k$, and such that
there exists $d<s<2d$ s.t.
\begin{equation*}
\liminf_{\|k\|_1\to\infty}\frac{P_k}{\|k\|_1^{-s}}>0.
\end{equation*}
Let $k_0$ be s.t. if $\|k\|_1>k_0$ then $P_k>0$, and let
\begin{equation*}
\gamma=\inf_{k>k_0}\frac{-\log(1-P_k)}{\|k\|_1^{-s}}>0.
\end{equation*}
Let $\rho>2(2k_0)^d$, and let $N$ and $l$ be integers. Let $C_1$ and $C_2$ be
cubes of side-length $N$, which are $l$ cubes away from each other.
Assume further that $C_1$ and $C_2$ contain clusters $U_1$ and $U_2$,
each of size $\rho N^\frac{s}{2}$.
Then, the probability that there is an open bond between a vertex in the
$U_1$ and a vertex in $U_2$ is at least $\eta_s(\zeta\gamma\rho^2,l)$
for $\zeta=2^{-s-1}d^{-s}$.
\end{lemma}


In order to prove Lemma \ref{normal1}, we will need a few definitions as well as another lemma, Lemma \ref{lem:aldous} below, which is proved
in Appendix \ref{app:aldous}.

Let $M$ be a (large) integer, and let $1<\xi<2$. An {\em inhomogenous random graph with size $M$ and parameter $\xi$}, as defined in \cite{aldous}, is a set of particles $h_1,\ldots,h_k$ of masses
$m(h_1),\ldots,m(h_k)$ such that $\sum_{i=1}^km(h_i)=M$, such that for every $i\neq j$, there is a bond between the particles $h_i$ and $h_j$ with probability
$\eta_\xi\big(m(h_i)\cdot m(h_j),M\big)$, and different bonds are independent of each other.

For a connected component $C$ in an inhomogenous random graph $H$, we say that its {\em mass} is $m(C)=\sum_{i:h_i\in C}m(h_i)$. For $\chi>0$ and an inhomogenous random graph $H$, we define
$N_\chi(H)$ to be the number of connected clusters in $H$ whose mass is greater than or equal to $\chi$. 

\begin{lemma}\label{lem:aldous}
Let $1<\xi<2$, and let $\gamma<1$ be such that 
\begin{equation}\label{eq:gammadef}
18\gamma>16+\xi.
\end{equation}
There exists $\varphi=\varphi(\xi,\gamma)>0$ and $M_0=M_0(\xi,\gamma)$ such that for all $M>M_0$ and
every inhomogenous random graph $H$ with size $M$ and parameter $\xi$,
\[
P\big( N_{M^\gamma}(H) \geq 2\big) < M^{-\varphi}.
\]
\end{lemma}

\begin{proof}[Proof of Lemma \ref{normal1}]
Let $s/d<\xi<2$, and let $\gamma<1$ be as in \eqref{eq:gammadef}.
Let $\varphi=\varphi(\xi,\gamma)$ be as in Lemma \ref{lem:aldous}.
Notice that by Theorem \ref {GKN} (\cite{uniq}) there is a unique infinite cluster.
Choose
\begin{equation*}
C_n=n^a\text{ and }D_n=n^{-b},
\end{equation*}
where 
\begin{equation}\label{eq:reqab}
a>b>1,\ 2b<a(2d-s), \mbox{ and }
(da-b)/da>\gamma.
\end{equation}
Choose $\epsilon'$ s.t.
\begin{equation}\label{prodconv}
3\epsilon'\prod_{k=1}^{\infty}{(1+3D_k)}<\epsilon.
\end{equation}
Such an $\epsilon'$ exists because the product in (\ref{prodconv}) converges.

Let
\begin{equation*}
\lambda = \inf_{\|k\|_1>0}\frac{-\log(1-P_k)}{\|k\|_1^{-s}}.
\end{equation*}
Notice that since
\begin{equation*}
\lim_{x\searrow 0}\frac{-\log(1-x)}{x}=1,
\end{equation*}
we get that $\lambda>0$. By the choice of $\lambda$,
for every $k$ 
we have that $P_k\geq\eta_s(\lambda,\|k\|_1)$.

Denote by $\alpha$ the density of the infinite percolation cluster. Let
$M>\max(M_0,2/\alpha)$ where $M_0$ is as in Lemma \ref{lem:aldous}
be s.t. the following conditions hold: 
\begin{enumerate}
\item
with probability bigger than $1-\epsilon '$ at
least $\frac{1}{2}\alpha M^d$
of the vertices in $[0,M-1]^d$ are in the infinite cluster.
\item
For every $n\geq 1$,
\begin{equation}\label{eq:Mgamma}
\frac \alpha 2M((n-1)!)^{da-b} \geq \big[M(n!)^{da}\big]^\gamma.
\end{equation}
\item
For every $n\geq 1$,
\begin{equation}\label{eq:Mxi}
\eta_\xi\left(1,2\big[M((n-1)!)^a\big]^d\right) \leq \eta_s(\lambda,dM(n!)^a)
\end{equation}
\item
For every $n\geq 1$,
\begin{equation}\label{eq:Mvarphi}
n^{2ad} < \big(M[(n-1)!]^a\big)^{\varphi/2}
\end{equation}
\item
\begin{equation}\label{eq:Mvarphieps}
M^{-d\varphi}<\epsilon'
\end{equation}
\item
For every $n$,
\begin{equation}\label{eq:Mdvarphieps}
\big[M((n-1)!^a)\big]^{-d\varphi/2}\leq\epsilon'D_n.
\end{equation}
\end{enumerate}
The existence of
this $M$ follows from the ($d$-dimensional) ergodic theorem as well as \eqref{eq:reqab} and the fact that
\eqref{eq:Mxi}, \eqref{eq:Mvarphi}, \eqref{eq:Mvarphieps} and \eqref{eq:Mdvarphieps} hold for all large enough $M$.
The infinite cluster is unique, so all of the percolating vertices in
$[0,M-1]^d$ will be connected to each other within some
big cube containing $[0,M-1]^d$.
Let $K$ be such that they are all connected inside $[-K,M+K-1]^d$ with
probability $>1-\epsilon '$.
We define a {\em semi-cluster} in a cube
\begin{equation*}
\C=\prod_{i=1}^{d}{[l_iM,(l_i+1)M-1]}\text{ }(l_i\in\Z\text{ }\forall i)
\end{equation*}
to be a maximal (w.r.t. containment) set of vertices in the cube that is contained in a connected subset of the
$K$-enlargement
\begin{equation*}
\C_K=\prod_{i=1}^{d}{[l_iM-K,(l_i+1)M+K-1]}.
\end{equation*}
of the cube.
We call a cube $\C$
{\em alive}
if there is a unique semi-cluster in $\C$ of size larger than $\frac{1}{2}\alpha M^d$.
We now show that by the choice of $M$ and $K$, the probability that there is more than one semi-cluster of size larger than $\frac{1}{2}\alpha M^d$ in
$\C$ is less than $\epsilon'$.
Indeed, for every $x,y\in\C$, we have $P_{x-y}>\eta_\xi(1,M^d)=:\upsilon_1$. Therefore, we can sample the configuration $\omega$ in two steps as follows: For every $x,y\in\C$, let $P'_{x,y}$ be the value such that $P'_{x,y}+\upsilon_1-\upsilon_1P'_{x,y}=P_{x-y}$. We then sample $\omega'$ as and independent configuration where the bond $(x,y)$ appears w.p. $P'_{x,y}$ if $x,y\in\C$ and $P_{x-y}$ otherwise.
We then sample the configuration $\omega''$ as an i.i.d. $\upsilon_1$ configuration on $\C$, and $\omega:=\omega'\cup\omega''$. Then $\omega$ has the required distribution (i.e. independent with probability $P_{x-y}$ for the edge $(x,y)$ for all $x$ and $y$).
Now, for every two $\omega'$-semi-clusters, $S_1$ and $S_2$, the probability that they are connected in $\omega$ is the probability that there is an $\omega''$ edge between them, which is $\eta_\xi(|S_1|\cdot|S_2|,M^d)$. Thus, the semi-clusters in $\omega'$ are an inhomogeneous random graph, and by Lemma \ref{lem:aldous}, the probability that there is more than one semi-cluster larger than $\frac{\alpha}{2}M^d$ is bounded by $M^{-d\varphi}<\epsilon'$.

Therefore, a cube of side-length $M$ is alive with probability at least $1-3\epsilon `$.


\ignore{
For every living cube, choose a semi-cluster (by {\em semi-cluster} we mean
a set of vertices in the cube that is contained in a connected subset of the
$K$-enlargement of the cube) of size at least $\frac{1}{2}\alpha M^d$ inside it.
We say that two cubes $\C_1$ and $\C_2$ are attached to each other if there exists an
open bond between the semi-cluster in $\C_1$ and the semi-cluster in $\C_2$.
If the cubes $\C_1$ and $\C_2$ are alive and are $k$ cubes
away from each other, then the probability that they are connected
is at least
\begin{equation*}
\eta_s(\gamma M^{2d-s},k)
\end{equation*}
for $\gamma=\frac{1}{4}\alpha^2\lambda(2d)^{-s}$.

This is true because there are at least $\frac{1}{4}\alpha^2M^{2d}$ pairs of vertices
$(v_1,v_2)$ from the semi-clusters of $C_1$ and $C_2$ respectively s.t. $\|v_1-v_2\|_1>k_0$.
For these vertices, $\|v_1-v_2\|_1<2dkM$. So, the probability that there is no edge between $v_1$
and $v_2$ is bounded by $1-\eta_s(\lambda,2dkm)=\exp(-\lambda(2dkM)^{-s})$. So, the probability
that there is no edge between the semi-cluster in $C_1$ and the one in $C_2$ is no more than
\begin{eqnarray*}
\left[\exp(-\lambda(2dkM)^{-s})\right]^{\frac{1}{4}\alpha^2M^{2d}}
&=&\exp\left(-\frac{1}{4}\alpha^2M^{2d}\lambda(2dkM)^{-s}\right)\\
&=&\exp\left(-\frac{1}{4}\alpha^2(2d)\lambda^{-s}M^{2d-s}k^{-s}\right)\\
&=&1-\eta_s(\gamma M^{2d-s},k).
\end{eqnarray*}
}

For $k=1,2,\ldots,$ let 
\[
M_k=M\prod_{l=1}^kC_l=M[k!]^a \ \ \ \ \ \mbox{and} \ \ \ \ \ U_k=\frac{\alpha}{2}\prod_{j=1}^{k-1}D_j=\frac{\alpha}{2}[(k-1)!]^{-b}
\]
Note that $M_1=M$.

\ignore{
Choose some large number $\beta$.
Take $M$
and $K$ s.t. $\gamma M^{2d-s}>\beta$ and s.t. the probability of a cube
to be alive is more than $1-2\epsilon'$. The probability that two cubes that are $k$
cubes away from each other are attached is at least $\eta_s(\beta,k)$.
}

Let $R$ be a number such that
\begin{equation}\label{katanmeps}
(MR+2K)^d<2(MR)^d.
\end{equation}
We want to renormalize to
cubes of side length $N=RM+K$.
We cannot apply the renormalization
argument from \cite{NS}, because the events that two (close enough) cubes
are alive are dependent. Thus, we use a different technique
of renormalization:

Consider the $M$-sided cubes as first stage vertices. Then, take cubes of side-length $C_1$ of
first stage vertices, and consider them as second stage vertices. Now, take cubes of side-length
$C_2$ of second stage vertices and consider them as third stage vertices.
Keep on taking cubes of
side length $C_n$ of $n$-stage vertices and consider them as $n+1$ stage vertices.

Choose $R$ to be
\begin{equation*}
R=\prod_{n=1}^{L}C_n
\end{equation*}
for $L$ large enough for (\ref{katanmeps}) to hold.

We already have a notion of a first stage vertex being alive.
Define inductively that an $n$-stage vertex is {\em alive} if
at least $D_n(C_n)^d$ of the $(n-1)$-stage vertices in it are alive, and every
two of those vertices are {\em attached} to each other, i.e.  there is an open bond
between the big clusters in these $n-1$ stage vertices.
Denote by $\lambda_n$ the
probability that an $n$-stage vertex is not alive. We want to bound $\lambda_n$:

Denote by $\phi_n$ the probability that there aren't enough living
$(n-1)$-stage vertices inside our $n$-stage vertex, and by $\psi_n$ the
probability that not every two of them are attached to each other. Then,
$\lambda_n\leq\phi_n+\psi_n$.
Given $\lambda_{n-1}$, the
expected number of dead $(n-1)$-stage vertices in
an $n$-stage vertex is $\lambda_{n-1}C_n^d$. Therefore, by the Markov
inequality,
\begin{equation*}
\phi_n\leq\frac{\lambda_{n-1}}{1-D_n}.
\end{equation*}

To estimate $\psi_n$, again we use the same argument using Lemma \ref{lem:aldous}. 
For every $\C_1$ and $\C_2$ of level $n-1$ vertices of distance bounded by $C_n$, and every $x,y\in\C=\C_1\cup\C_2$, we have $P_{x-y}>\eta_\xi(1,2M_{n-1}^d)=:\upsilon_{n}$. Therefore, we can sample the configuration $\omega$ in two steps as follows: For every $x,y\in\C$, let $P'_{x,y}$ be the value such that $P'_{x,y}+\upsilon_n-\upsilon_nP'_{x,y}=P_{x-y}$. We then sample $\omega'$ as and independent configuration where the bond $(x,y)$ appears w.p. $P'_{x,y}$ if $x,y\in\C$ and $P_{x-y}$ otherwise.
We then sample the configuration $\omega''$ and and i.i.d. $\upsilon_n$ configuration on $\C$, and $\omega:=\omega'\cup\omega''$. Then $\omega$ has the required distribution (i.e. independent with probability $P_{x-y}$ for the edge $(x,y)$ for all $x$ and $y$).
Now, for every two $\omega'$-semi-clusters, $S_1$ and $S_2$, the probability that they are connected in $\omega$ is the probability that there is an $\omega''$ edge between them, which is $\eta_\xi(|S_1|\cdot|S_2|,M^d)$. Thus, the semi-clusters in $\omega'$ are an inhomogeneous random graph, and by Lemma \ref{lem:aldous}, the probability that there is more than one semi-cluster larger than $U_nM_{n-1}^d$ is bounded by $M_{n-1}^{-d\varphi}$. However, if both $\C_1$ and $\C_2$ are alive and the big components in them are not connected to each other, then there are (at least) two semi-clusters in $\C$ which are larger than $U_nM_{n-1}^d$. Therefore,
\begin{eqnarray*}
\psi_n \leq P\left[
\exists_{\C_1 \mbox{ and } \C_2} \mbox{in $[0,M_n]^d$, alive and not connected}
\right]
\leq M_{n-1}^{-d\varphi}{{M_n^d}\choose{2}}\leq M_{n-1}^{-d\varphi/2}
\leq\epsilon'D_n.
\end{eqnarray*}

\ignore{
Every living $(n-1)$-stage vertex includes at least
\begin{equation*}
V_n=\prod_{k=1}^{n-1}(C_k)^dD_k=((n-1)!)^{da-b}
\end{equation*}
living first-stage vertices inside its connected component. The distance
between those first-stage vertices cannot exceed
\begin{equation*}
U_n=d\prod_{k=1}^{n}C_k=d(n!)^a.
\end{equation*}
Therefore,
\begin{equation*}
\psi_n\leq \binom{(C_n)^d}{2}(1-\eta_s(\beta,U_n))^{{V_n}^2}
\leq {C_n}^{2d}e^{-\beta U_n^{-s}{V_n}^2}
\end{equation*}
i.e.
\begin{eqnarray*}
\psi_n\leq n^{2da}e^{-\beta(d^{-s}(n!)^{-as}\cdot ((n-1)!)^{2(da-b)})}\\
=n^{2ad}\cdot e^{-\beta(d^{-s}n^{-as}\cdot ((n-1)!)^{a(2d-s)-2b})}.
\end{eqnarray*}
(Notice that the event that the connecting edges exist might depend on the
existence of enough living vertices. However, in this case, the FKG inequality
works in our favor)
\\This shows that $\psi_n$ decays faster than exponentially, and therefore,
since we control $\beta$ and can make it as large as we like, we can achieve
\begin{equation*}
\psi_n<\epsilon 'D_n
\end{equation*}
for every $n$.
}

By the choice of $M$ and $K$, and by the definition of $\lambda_1$, we see
that $\lambda_1<3\epsilon '$. In addition, for every $n$,
\begin{eqnarray*}
\lambda_n\leq\psi_n+\phi_n
&\leq&\epsilon 'D_n+\frac{\lambda_{n-1}}{1-D_n}\\
&\leq&\epsilon 'D_n+\lambda_{n-1}(1+2D_n)\\
&\leq&(1+3D_n)\max(\lambda_{n-1},\epsilon')
\end{eqnarray*}
Therefore, by induction, we get that for every $n$
\begin{equation*}
\lambda_n\leq 3\epsilon '\prod_{k=1}^{n}(1+3D_k),
\end{equation*}
and so, for all $n$,
\begin{equation*}
\lambda_n\leq \Theta\epsilon '
\end{equation*}
 where
\begin{equation*}
\Theta=2\prod_{k=1}^{\infty}{(1+3D_k)}<\infty.
\end{equation*}

So, with probability at least $1-\Theta\epsilon '>1-\epsilon$, we have a
cluster of size
\begin{equation*}
\prod_{n=1}^L{D_n(C_n)^d}=\prod_{n=1}^{L}{n^{da-b}}
=\left(\prod_{n=1}^L{C_n}\right)^\frac{da-b}{a}
=R^\frac{da-b}{a}.
\end{equation*}

This is larger than $2\rho R^{\frac{s}{2}}$
If $L$ is large enough, because
\begin{equation*}
\frac{da-b}{a}>\frac{s}{2}.
\end{equation*}

So, by (\ref{katanmeps}), the lemma is proved for
$N=RM+2K$.
\end{proof}

\begin{proof}[Proof of Lemma \ref{connectprob} ]
There are $\rho^2N^s$ pairs of vertices $(v_1,v_2)$ s.t. $v_1\in U_1$ and $v_2\in U_2$.
For every $v_1\in U_1$ there are at most $(2k_0)^d<\halb\rho N^\frac{s}{2}$ vertices at
distance smaller or equal to $k_0$ from $v_1$. So, at least half of the pairs $(v_1,v_2)$
satisfy $\|v_1-v_2\|_1>k_0$. All of the pairs satisfy $\|v_1-v_2\|_1\leq 2ldN$.
For a given pair $(v_1, v_2)$ s.t. $\|v_1-v_2\|_1>k_0$,
the probability that there is no edge between $v_1$ and $v_2$ is
bounded by $1-\eta_s(\gamma,2ldN)$. So, the probability that
there is no edge between $U_1$ and $U_2$ is bounded by
\begin{eqnarray*}
\left[1-\eta_s(\gamma,2ldN)\right]^{\halb\rho^2N^s}
&=&\left[\exp(-\gamma(2ldN)^{-s})\right]^{\halb\rho^2N^s}\\
&=&\exp(-\gamma(2ldN)^{-s}\cdot\halb\rho^2N^s)\\
&=&\exp(-\halb(2d)^{-s}\gamma\rho^2l^{-s})\\
&=&1-\eta_s(\zeta\gamma\rho^2,l)
\end{eqnarray*}
\end{proof}

We can now use Lemma \ref{normal1} and Lemma \ref{connectprob} to prove the following
extension of Theorem \ref{criti}:
\begin{theorem}
Let $d\geq 1$, and let $\{P_k\}_{k\in\Z^d}$ be probabilities such that there exists $s<2d$
for which
\begin{equation}\label{hasum}
\liminf_{\|k\|\to\infty}\frac{P_k}{{\|k\|_1^{-s}}}>0.
\end{equation}
Then, if $\{P_k\}$ is percolating then there exists an $\epsilon>0$ such that
$\{P'_k=(1-\epsilon)P_k\}$ is percolating too.
\end{theorem}


\begin{proof}
Let $\{P_k\}_{k\in\Z^d}$ be a percolating system that satisfies (\ref{hasum}).
Let $k_0$ and $\gamma$ be as in Lemma \ref{connectprob}. Let, again,
$\zeta=2^{-s-1}d^{-s}$.

Let $\lambda<1$, $\beta$ and $\delta>0$ be such that a system in which every
vertex is alive with probability $\lambda-\delta$ and every two vertices $x$
and $y$ are connected to each other with probability $\eta_s(\beta(1-\delta),\|x-y\|_1)$
is percolating. For one dimension one can choose such $\lambda, \beta$ and $\delta$
by Theorem \ref{thm:ns}. For higher dimensions we may use the fact that site-bond
nearest neighbor percolation with high enough parameters has, a.s., an infinite cluster.

Let $\rho>2(2k_0)^d$ be s.t. $\zeta\gamma\rho^2\geq\beta$. By Lemma \ref{normal1},
there exists $N$ s.t.
a cube of side length $N$ contains a cluster of size $\rho N^\frac{s}{2}$
with probability bigger than $\lambda$. A Cube that contains
a cluster of size bigger or equal to $\rho N^\frac{s}{2}$ will be considered {\em alive}.
For $\omega\leq 1$, consider the system $\{P'_k=\omega P_k\}$.
The probability that in the system $\{P'_k\}$
an $N$-cube is alive is a continuous function of $\omega$.
If we define $k'_0$ and $\gamma'$ for $\{P'_k\}$
the same way we defined $k_0$ and $\gamma$, then we get that $k'_0=k_0$, and
$\gamma'$ is a continuous function of $\omega$.

Choose $\epsilon$ be so small that in the system $\{P'_k=(1-\epsilon)P_k\}$ the
probability of an $N$-cube to be alive is no less than $\lambda-\delta$ and that
$\gamma'\geq(1-\delta)\gamma$. Then, in the system
$\{P'_k\}$, every $N$-cube is alive with probability bigger than $\lambda-\delta$, and
two cubes at distance $k$ cubes from each other are connected with probability
bigger than
\begin{eqnarray*}
\eta_s(\zeta\gamma'\rho^2,k)&=&\eta_s((1-\delta)\zeta\gamma\rho^2,k)\\
&\geq&\eta_s(\beta(1-\delta),k).
\end{eqnarray*}
So, by the choice of $\beta$, $\lambda$ and $\delta$, a.s. there is an infinite cluster
in the system $\{P'_k\}$.
\end{proof}

\begin{corollary}
Let $d\geq 1$, and let $\{P_k\}_{k\in\Z^d}$ be probabilities such that there exists $s<2d$
for which
\begin{equation}
\liminf_{\|k\|\to\infty}\frac{P_k}{{\|k\|_1^{-s}}}>0.
\end{equation}
If $\{P_k\}$ is critical, i.e. for every $\epsilon>0$ the system
$\{(1+\epsilon)P_k\}$
is percolating but the system
$\{(1-\epsilon)P_k\}$
is not percolating, then $\{P_k\}$ is not percolating.
\end{corollary}

Lemma \ref{normal1} also serves us in proving Theorem \ref{jeff}.
\begin{proof}[Proof of Theorem \ref{jeff}]
Let $\{P_k\}_{k\in Z^d}$ be such that
\begin{equation*}
\liminf_{\|k\|\to\infty}\frac{P_k}{\|k\|_1^{-s}}>0
\end{equation*}
for $s<2d$.
Let $k_0$, $\gamma$ and $\zeta$ be as before.
Let $\epsilon$ and $\rho>2(2k_0)^d$ be s.t. site-bond nearest neighbor
percolation s.t. every site is alive with probability $1-\epsilon$
and every bond is open with probability $\eta_s(\zeta\gamma\rho^2,1)$ on $Z^d$
percolates. Let $N$ be suitable for those $\epsilon$ and $\rho$ by
Lemma \ref{normal1}. Now, erase all of the bonds of length bigger than $4Nd$.
Renormalize the space to cubes of side-length $N$.
By erasing only bonds of length $>4Nd$, we did not erase bonds that are
either inside
$N$-cubes, or between neighboring $N$-cubes.
So, the renormalized picture still gives us site-bond percolation with probabilities
$1-\epsilon$ and $\eta_s(\zeta\gamma\rho^2,1)$, and therefore an infinite cluster exists a.s.
\end{proof}

Returning to transience, we now prove that for large enough parameters $\beta$
and $\lambda$, the infinite cluster is transient. Later we will use
Lemma \ref{normal1} and Lemma \ref{connectprob} to reduce any percolating system
(with $d<s<2d$) to one with these large $\beta$ and $\lambda$.
\begin{lemma}\label{trans1}
Let $d\geq 1$ and $d<s<2d$. Consider the independent bond-site percolation model in which every
two vertices, $i$ and $j$, are connected with probability
$\eta_s(\beta,\|i-j\|_1)$, and every vertex is alive with probability $\lambda<1$.
If $\beta$ is large enough and $\lambda$ is close enough to $1$, then (a.s.) the random
walk on the infinite cluster is transient.
\end{lemma}

In order to prove Lemma \ref{trans1}, we need the notion of a {\em renormalized graph}:
For a sequence $\{C_n\}_{n=1}^{\infty}$, we construct a graph whose vertices are marked
$V_l(j_l,..,j_1)$ where $l=0,1,...$ and $1\leq j_n\leq C_n$. For convenience,
set $V_k(0,0,..,0,j_l,..,j_1)=V_l(j_l,..,j_1)$. For $l\geq m$,
we define $V_l(j_l,...,j_m)$ to be the set
\begin{equation*}
\{V_l(j_l,...,j_m,u_{m-1},...,u_1)|1\leq u_{m-1}\leq C_{m-1},...,1\leq u_1\leq C_1\}.
\end{equation*}

\begin{definition}
A {\em renormalized graph} for a sequence $\{C_n\}_{n=1}^{\infty}$ is a graph whose vertices
are $V_l(j_l,..,j_1)$ where $l=0,1,...$ and $1\leq j_n\leq C_n$, such that
for every $k\geq l>2$, every $j_k,...,j_{l+1}$ and every $u_l,u_{l-1}$ and
$w_l,w_{l-1}$, there is an edge connecting a vertex in $V_k(j_k,...,j_{l+1},u_l,u_{l-1})$
and a vertex in $V_k(j_k,...,j_{l+1},w_l,w_{l-1})$.
\end{definition}


one may view a renormalized graph as a graph having the following recursive structure:
The $n$-th stage of the graph is composed of $C_n$
graphs of stage $(n-1)$, such that every $(n-2)$-stage graph in each of
them is connected to every $(n-2)$-stage graph in any other.
(A zero stage graph is a vertex).

\begin{lemma}\label{subgraph}
Under the conditions of Lemma \ref{trans1}, if $\beta$
and $\lambda$ are large enough, then a.s the infinite cluster contains a
renormalized sub-graph with $C_n=(n+1)^{2d}$.
\end{lemma}

\begin{proof}
We will show
that with a positive probability $0$ belongs to a renormalized sub-graph. Then, by
ergodicity of the shift operator and the fact that the event
$E=\{\text{There exists a renormalized sub-graph}\}$ is shift invariant we get
$\prob(E)=1$. In order to do that, we use the exact same technique used by Newman
and Schulman in \cite{NS}:

Take
\begin{equation}\label{hagda}
W_n=2(n+1)^2 \re ; \re\theta_n=1-\frac{n^{-1.5}}{2} \re;\re \lambda_n=1-\frac{(n+1)^{-1.5}}{4}.
\end{equation}
Renormalize $Z^d$ by viewing cubes of side-length $W_1$ as {\em first stage} vertices.
(the original vertices will be viewed as zero-stage vertices). Then, take cubes
of side-length $W_2$ of first-stage vertices as {\em second stage} vertices,
and continue grouping together cubes of side-length $W_n$ of $(n-1)$-stage vertices
to form {\em $n$ stage} vertices.

We now define inductively the notion of an ($n$-stage) vertex being alive:
The notion of a zero-stage vertex being alive is given to us. A first-stage
vertex is {\em alive} if at least $\theta_1W_1^d$ of its vertices are alive,
and they are all connected to each other.
For every living first-stage vertex, we choose $C_1$ zero-stage vertices, and call
them {\em active}. The {\em active part} of a first-stage vertex is the set of active
zero-stage vertices in it.
The active part of a living zero-stage vertex is the singleton containing the vertex.

We now define (inductively) simultaneously the notion of an $n$-stage vertex being alive,
and of the active part of this vertex.

For $n\geq 2$, we say that an $n$-stage vertex $v$ is {\em alive} if:
\\(A) At least $\theta_nW_n^d$ of its vertices are alive, and
\\(B) If $i_1$ is a living $(n-2)$-stage vertex that belongs to a living
$(n-1)$-stage vertex $i_2$ that belongs to $v$, and $j_1$ is a living $(n-2)$-stage
vertex that belongs to a living $(n-1)$-stage vertex $j_2$ that belongs to $v$ then
there exists an open bond connecting a zero-stage vertex in the active part of $i_1$
to a zero-stage vertex in the active part of $j_1$.

When choosing the active vertices, if the vertex that includes $0$ is alive, we choose it to
be active.

To define the active part: If $v$ is a living $n$-stage vertex, then we choose $C_n$ of its
living $n-1$-stage vertices to be active.
The active part of $v$ is the union of the active parts
of its active vertices. (notice that the active part is always a set of zero-stage vertices).

We denote the event that (A) occurs for the $n$-stage vertex containing $0$ by $A_n$,
and by $B_n$ we denote that (B) occurs for the $n$-stage vertex containing $0$.
$A_n(v)$ and $B_n(v)$ will denote the same event for the $n$-stage vertex $v$. Of
course, \prob($A_n)=\prob(A_n(v))$ and $\prob(A_n)=\prob(A_n(v))$ for every $v$.
Further, we denote by $L_n(v)$ the event that the $n$-stage vertex $v$ is alive,
and by $L_n$ the event that the $n$-stage vertex containing $0$ is alive.

Let $v$ be an $n$-stage vertex. Given $A_n$ we want
to estimate the probability of $B_n$: We have
at most
\begin{equation}\label{num_pairs}
{{(W_nW_{n-1})^d}\choose{2}}<4^d(n+1)^{8d}
\end{equation}
pairs of $(n-2)$-vertices.

If $i_1$ and $i_2$ are living $(n-2)$-stage vertices in $v$, then the distance between
a zero-stage vertex in $i_1$ and a zero-stage vertex in $i_2$ cannot exceed
\begin{equation}\label{dist_act}
\prod_{k=1}^{n}{W_k}=2^n((n+1)!)^2.
\end{equation}

The size of the active part in $i_1$ (and in $i_2$), is
\begin{equation}\label{num_act}
\prod_{k=1}^{n}{W_k}=((n+1)!)^{2d}
\end{equation}

By (\ref{num_act}) and (\ref{dist_act}), the probability that there is no open bond
between $i_1$ and $i_2$ is bounded by
\begin{eqnarray*}
\left[\exp\left(-\beta\cdot 2^{-ns}((n+1)!)^{-2s}\right)\right]^{((n+1)!)^{4d}}\\
=\exp\left(-\beta\cdot 2^{-ns}((n+1)!)^{4d-2s}\right)
\end{eqnarray*}
and by (\ref{num_pairs}) we get
\begin{eqnarray}\label{probB}
\prob\left[B_n^c|A_n\right]
 &\leq&
4^d(n+1)^{8d}\exp\left(-\beta\cdot 2^{-ns}((n+1)!)^{4d-2s}\right)\\
\nonumber &\leq& \exp\left(9d\log(n)-\beta\cdot 2^{-ns}((n+1)!)^{4d-2s}\right)
\end{eqnarray}
Assume that $\beta>1$. We may assume that because we deal with
"large enough" $\beta$.
By (\ref{probB}), there exists $n_0$ s.t. if $n>n_0$ then
\begin{equation}\label{fastthanexp}
\prob\left[B_n^c|A_n\right]<e^{-n}.
\end{equation}

We now want to prove the following claim:
\begin{claim}\label{induc}
There exists $n_1$ such that for every $n>n_1$, if $\prob(L_n)\geq\lambda_n$ then
$\prob(L_{n+1})\geq\lambda_{n+1}$.
\end{claim}
\begin{proof}
Let $\psi=\prob(L_n)$. First, we like to estimate $\prob(A_{n+1})$. The event
$A_{n+1}^c$ is the event that at least $(1-\theta_{n+1})W_{n+1}^d$ vertices are
dead. The number of dead vertices is a $(W_{n+1}^d, \psi)$ binomial variable,
and by the induction hypothesis together with (\ref{hagda}),
$\psi < \halb(1-\theta_{n+1})$. So, by large deviation estimates,
\begin{eqnarray}\label{larged}
\prob(A_{n+1}^c) &<& \exp\left(-\frac{\halb(1-\theta_{n+1})W_{n+1}^d}{16}\right)\\
\nonumber &\leq& \exp\left(-\frac{n^{2d-1.5}}{32}\right)
\end{eqnarray}
If $n_1>n_0$, and is large enough, by (\ref{fastthanexp}) and (\ref{larged}),
\begin{eqnarray*}
\prob(L_{n+1}^c) &\leq& \prob(A_{n+1}^c)+\prob(B_{n+1}^c|A_{n+1})\\
&\leq& \exp\left(-\frac{n^{2d-1.5}}{32}\right)+e^{-n}\\
&\leq& \frac{(n+1)^{-1.5}}{4}=1-\lambda_{n+1}
\end{eqnarray*}

\end{proof}

We can take $\beta$ and $\lambda$ so large that $\prob(L_{n_1})>\lambda_{n_1}$. But then,
by Claim \ref{induc}, for every $n>n_1$, $\prob(L_{n})>\lambda_{n}$. So, since the events
$L_n$ are positively correlated,
\begin{equation*}
\prob\left(\bigcap_{n=1}^{\infty}L_n\right)\geq\prod_{n=1}^{\infty}{\prob(L_n)}>0
\end{equation*}
So with positive probability, $0$ is in an infinite cluster. The active part of the infinite
cluster (i.e. the union of the active parts of the $n$-stage vertex containing $0$ for all $n$)
is a renormalized sub-graph of the infinite cluster that contains $0$.
\end{proof}

\begin{proof}[proof of Lemma \ref{trans1}]
In view of Lemma \ref{subgraph} it suffices to show that for $C_n=(n+1)^{2d}$,
the renormalized graph is transient.

We build, inductively, a flow $F$ from $V_1(1)$ to infinity which has a
finite energy. First, $F$ flows $C_1^{-1}$ mass from $V_1(1)$ to each of
\begin{equation*}
\{V_1(i)\}_{i=2}^{C_1}.
\end{equation*}
Now, inductively, assume that $F$ distributes the mass among
\begin{equation*}
\{V_n(i_1,...,i_n)|2\leq i_k\leq C_k\}.
\end{equation*}
Then, for each $(n-1)$-stage graph $V_n(i),i\neq 1$ and every $n$-stage graph
$V_{n+1}(j),j\neq 1$, there are two vertices, $p^{(n)}_{i,j}\in V_n(i)$ and
$q^{(n)}_{i,j}\in V_{n+1}(j)$ which are connected to each other by an open bond.
(Notice that the vertices $\{p^{(n)}_{i,2},...p^{(n)}_{i,C_{n+1}}\}$, as well as
$\{q^{(n)}_{2,j},...q^{(n)}_{C_n,j}\}$ do not necessarily differ from each other).
Inductively, we know how to flow mass from one vertex in $V_n(i)$ to all of
$V_n(i)$. We can flow it backwards in the same manner to any desired vertex.
Flow the mass so that it will be distributed equally among
$\{p^{(n)}_{i,2},...p^{(n)}_{i,C_{n+1}}\}$ (if a vertex appears twice, it will get a
double portion). Now flow the mass from each $p^{(n)}_{i,j}$ to the corresponding
$q^{(n)}_{i,j}$,
and from $q^{(n)}_{i,j}$ (again by the inductive familiar way) we will flood
$V_{n+1}(j)$. Now, we bound the energy of the flow: $E_n$, The maximal
possible energy of the first $n$ stages of the flow (i.e. the part of the flow
which distributes the mass the origin to $V_{n+1}$ and takes it backwards to
$\{p^{(n+1)}_{i,j}\}\subset V_{n+1}$) can be
bounded by the energy of first $n-1$ stages of the flow, plus:
\\(A) Flowing between $p^{(n)}_{i,j}$ and $q^{(n)}_{j,i}$: This will have energy of
$(C_nC_{n+1})^{-1}$.
\\(B) Flowing inside $V_{n+1}$: the energy is bounded by
$\frac{E_{n-1}}{C_{n+1}}$.
\\So,
\begin{equation*}
E_n\leq \left(1+\frac{1}{C_{n-1}}\right)E_{n-1}+\frac{1}{C_nC_{n-1}}.
\end{equation*}
The total energy is bounded by the supremum of $\{E_n\}$ which is finite because
\begin{equation*}
\sum_{n=1}^{\infty}{\frac{1}{C_n}}<\infty.
\end{equation*}
\end{proof}

Let $v$ be a vertex. The amount of flow that goes through $v$ is defined to be
$
f(v)=\halb\sum{|f(e)|}
$
where the sum is taken over all of the edges $e$ that have $v$ as an end point.
Then, we get a notion of the {\em energy
of the flow through the vertices}, defined as

\begin{equation*}
\en_{\text{vertices}}=\suml_{v\text{ is a vertex}}f(v)^2
\end{equation*}

\begin{remark}\label{vertices}
The same calculation as in Lemma \ref{trans1} yields that not only the energy of the flow on the
bonds is finite, but also the energy of the flow through the vertices.
\end{remark}

This fact allows us to obtain the the main goal of this section:

\begin{theorem}\label{thm:transience}
Let $d\geq 1$, and let $\{P_k\}_{k\in \Z^d}^\infty$ satisfy:
\\(A) $P_k = P_{-k}$ for every $k\in\Z$.
\\(B) the independent percolation model in which the bond between i and j is
open with probability $P_{i-j}$ has, a.s., an infinite cluster.
\\(C) there exists $d<s<2d$ s.t.
\begin{equation}\label{katanm2}
\liminf_{\|k\|\to\infty}\frac{P_k}{{\|k\|_1^{-s}}}>0.
\end{equation}
\\(D)
\begin{equation*}\label{finitdeg}
\sum_{k\in\Z^d}{P_k}<\infty.
\end{equation*}

Then, a.s., a random walk on the infinite cluster is transient.
\end{theorem}
\begin{proof}
By (D), the degree of every vertex in the infinite cluster is finite, so the random walk
is well defined.

Let $\beta$ and $\lambda$ be
large enough for Lemma \ref{trans1}. Then, by Lemma \ref{normal1}, there
exists $N$ such that after renormalizing with cubes of side-length $N$ we
get a system whose connection probabilities dominate
$\eta_s(\beta,|i-j|)$, and the probability of a vertex to live is bigger than
$\lambda$. By Lemma \ref{trans1}, there is a flow on this graph whose energy is
finite. For the walk to be transient, the energy of the flow should also be
finite inside the $N$-cubes. This is true because of Remark
\ref{vertices} and the fact that inside each $N$-cube there are no
more than
$\left(
\begin{smallmatrix}
N^d
\\2
\end{smallmatrix}
\right)$
bonds.
\end{proof}

one can look on other types of energy as well. For any $q$, we define the
$q$-energy of a flow as in equation (\ref{energy}).
Theorem \ref{thm:transience} says that for every $\{P_k\}$ that satisfies conditions
(A) through (D),
there is a flow
with finite $2$-energy. Actually, one can say more:
\begin{theorem}
Let $\{P_k\}_{k\in\Z}$ be as in Theorem \ref{thm:transience}. Then, For
every $q>1$, there is a flow with finite $q$-energy on the infinite cluster.
\end{theorem}

\begin{proof}[A sketch of the proof]
The proof is essentially the same as the proof of Theorem \ref{thm:transience}.
We can construct a renormalized sub-graph of the infinite cluster with
$C_n=(n+1)^{kd}$, for $k$ s.t. $k(q-1)>1$. We construct the flow the same way we did
it in Lemma \ref{trans1}. The same energy estimation will now yield the
required finiteness of the energy. Lemma \ref{normal1} and
Remark \ref{vertices} are used the same way they were used in Theorem
\ref{thm:transience}.

If we construct a renormalized graph with $C_n=2^n$ (such a graph a.s. exists as a
sub-graph
of the infinite cluster), we get a flow whose $q$-energy is finite for every $q>1$.
\end{proof}


\section{The recurrence proofs}\label{recprf}

In this section we prove the recurrence results. Unlike the transient case,
here we give two different proofs - one for the one-dimensional case, and the
other for the two-dimensional case. We begin with the easier
one-dimensional case.

\begin{theorem}\label{one_d_recur}
Let $\{P_k\}_{k=1}^\infty$ be a sequence of probabilities s.t.:
\\(A) the independent percolation model in which the bond between i and j is
open with probability $P_{|i-j|}$ has, a.s., an infinite cluster, and
\\(B)
\begin{equation*}
\limsup_{k\to\infty}\frac{P_k}{k^{-2}}<\infty.
\end{equation*}
Then, a.s., a random walk on the infinite cluster is recurrent.
\end{theorem}
The proof of the theorem relies on the Nash-Williams theorem, whose proof
can be found in \cite{yuval}:
\begin{theorem}[Nash-Williams]\label{Thm:nw}
Let $G$ be a graph with conductance $C_e$ on every edge $e$. Consider a
random walk on the graph such that when the particle is at some vertex, it
chooses its way with probabilities proportional to the conductances on the
edges that it sees.
Let $\{\Pi_n\}_{n=1}^\infty$ be disjoint cut-sets, and Denote by
$C_{\Pi_n}$ the sum of the conductances in $\Pi_n$.
If
\begin{equation*}
\sum_n{C_{\Pi_n}^{-1}}=\infty
\end{equation*}
then the random walk is recurrent.
\end{theorem}
In order to prove theorem \ref{one_d_recur}, we need the following definition and
three easy lemmas.

The following definition appeared originally in \cite{AN} and \cite{NS}.
\begin{definition}[Continuum Bond Model]\label{defcont}
Let $\beta$ be s.t.
\begin{equation*}
\integral_{0}^{1}\integral_{k}^{k+1}{\beta(x-y)^{-2}dydx}>P_k
\end{equation*}
for every $k$. The {\em continuum bond model} is the two dimensional inhomogeneous Poisson
process $\xi$ with density $\beta(x-y)^{-2}$. We say that two sets $A$ and $B$ are
connected if $\xi(A\times B)>0$.
\end{definition}
Notice that by the definition \ref{defcont}, the probability that the interval $[i,i+1]$ is
connected to $[j,j+1]$ in the continuum model is not smaller than the probability that $i$
is directly connected to $j$ in the original model. (By saying that a vertex is
{\em directly connected} to an interval , we mean
that there is an open bond between this vertex and some vertex in the interval.) So, we get:
\begin{claim}\label{condisc}
Let $I$ be an interval. Let $M$ be the length of the shortest interval that contains all of
the vertices that are directly connected to $I$ in the original model. Let $M'$ be the length
of the smallest interval $J$ s.t. $\xi(I\times (\R-J))=0$. Then, $M'$ stochastically dominates
$M$.
\end{claim}

\begin{lemma}\label{soi}
(A) Under the conditions of Theorem \ref{one_d_recur}, let $I$ be an interval of
length $N$. Then, the probability that there exists a vertex of distance
bigger than $d$ from the interval, that is directly connected to the interval, is
$O\left(\frac{N}{d}\right)$.\\
(B) Consider the continuum bond model. Let $I$ be an interval of length $N$, and let
$J$ be the smallest interval s.t. $\xi(I\times (\R-J))=0$. Then
$\prob(|J|>d)=O\left(\frac{N}{d}\right)$.
\end{lemma}

\begin{proof}
(A) Let
\begin{equation*}
\beta'=\sup_{k}\frac{P_k}{k^{-2}}<\infty.
\end{equation*}
If $v$ is at distance $k$ from $I$, then the probability that $d$ is directly
connected to $I$ is bounded by
\begin{equation*}
\beta'\sum_{k=d}^{d+N}{k^{-2}}<\frac{\beta' N}{d^2}
\end{equation*}
So, the probability that there is a vertex of distance bigger than $d$ that is directly
connected to $I$ is bounded by
\begin{equation*}
\sum_{k=d}^{\infty}\frac{\beta' N}{k^2} = O\left(\frac{N}{d}\right)
\end{equation*}
(B) is proved exactly the same way.
\end{proof}

\begin{lemma}\label{log}
Under the same conditions, and again letting $I$ be an interval of length $N$,
the expected number of open bonds exiting $I$ is O($\log N$).
There is a constant $\gamma$, s.t. the probability of having
more than $\gamma\log N$ open bonds exiting $I$ is smaller than $0.5$.
\end{lemma}

\begin{proof}
Again, let
\begin{equation*}
\beta'=\sup_{k}\frac{P_k}{k^{-2}}<\infty.
\end{equation*}
The expected number of open bonds exiting $I$ is
\begin{eqnarray*}
\sum_{v\in I, u\notin I}\prob(v\leftrightarrow u)&\leq&
\beta'\sum_{v\in I, u\notin I}(u-v)^{-2}\\
&=&2\beta'\sum_{i=1}^{N}\sum_{k=i}^{\infty}k^{-2}\\
&\leq&4\beta'\sum_{i=1}^{N}\frac{1}{i}\\
&=&O(\log N).
\end{eqnarray*}
Let $C$ be s.t. the expected value is less than $C\log N$ for all $n$. For any $\gamma>2C$,
by Markov's inequality, the probability that more than $\gamma\log N$ open bonds are
exiting $I$ is smaller than $0.5$.
\end{proof}

\begin{lemma}\label{halfovern}
Let $A_i$ be independent events s.t. $\prob(A_i)\geq 0.5$ for every $i$.
Then, a.s.,
\begin{equation*}
\sum_{i=1}^{\infty}\frac{1_{A_n}}{n}=\infty
\end{equation*}
\end{lemma}

\begin{proof}
Let
\begin{equation*}
U_k=\sum_{i=2^k}^{2^{k+1}-1}\frac{1_{A_i}}{i}
\end{equation*}
Then,
\begin{equation}\label{U}
U_k\geq 2^{-(k+1)}\sum_{i=2^k}^{2^{k+1}-1}1_{A_i}
\end{equation}
The variables $U_k$ are independent of each other, and by (\ref{U}),
for every $k$ we have $\prob(U_k>0.25)>0.5$. Therefore,
\begin{equation*}
\sum_{n=1}^{\infty}\frac{1_{A_n}}{n}=\sum_{k=0}^{\infty}U_k=\infty
\end{equation*}
a.s.
\end{proof}

\begin{proof}[Proof of theorem \ref{one_d_recur}]
We will show that with probability $1$, the infinite cluster satisfies the
Nash-Williams condition.
Let $I_0$ be some interval. We define $I_n$ inductively to be the smallest
interval that contains all of the vertices that are connected directly to $I_{n-1}$.
Denote
\begin{equation*}
D_n=\frac{|I_{n+1}|}{|I_n|}.
\end{equation*}
The edges exiting $I_{n+1}$ are stochastically dominated by the edges exiting
an interval of length $|I_{n+1}|$. (without the restriction that no edge starting
at $I_n$ exits $I_{n+1}$). Furthermore, given $I_n$ the edges exiting
$I_{n+1}$ are independent of those exiting $I_n$.
Let $\{U_n\}_{n=1}^\infty$ be independent copies of the continuum bond model.
Then, by Claim \ref{condisc} $D_n$ is stochastically dominated by the sequence
$D'_n=\frac{|I'_{n+1}|}{|I_n|}$,
where $I'_{n+1}$ is the smallest interval s.t. $\R-I'_{n-1}$ is not connected to
the copy of $I_{n}$ in $U_n$.

The variables $D'_n$ are i.i.d. Therefore, by Lemma \ref{soi}, the sequence
$\{\log(D_n)\}$ is dominated by a sequence of
i.i.d. variables $d_n=\log(D'_n)$, which satisfy $\E (d_n)<M$.
Let $\Pi_n$ be the set of
bonds exiting $I_n$. Then, $\{\Pi_n\}_{n=1}^\infty$ are disjoint cut-sets.
Given the intervals $\{I_n\}_{n=1}^N$, the set $\Pi_N$ is independent of
$\{\Pi_n\}_{n=1}^{N-1}$. Now, independently for each $n$, by Lemma \ref{log},
with probability bigger than $0.5$,
\begin{equation}\label{pin}
|\Pi_N|<\gamma\sum_{n=1}^{N}{d_n}.
\end{equation}
By the strong law of large numbers, with probability $1$, for all large
enough $N$,
\begin{equation}\label{sumdn}
\sum_{n=1}^{N}{d_n}<2MN.
\end{equation}
Combining (\ref{pin}), (\ref{sumdn}) and Lemma \ref{halfovern}, we get that
the Nash-Williams condition is a.s. satisfied.
\end{proof}

We now work on the two-dimensional case. Our strategy in this case will
be to project the long bonds on the short ones. That is, for every open long bond we
find a path of nearest-neighbor bonds s.t. the end points of the path are those of
the original long bond. Then, we erase the long bond, and assign its conductance to
this path. In order to keep the conductance of the whole graph, if the path is
of length $n$, we add $n$ to the conductance of each of the bonds involved
in it. To make the discussion above more precise, we state it as a lemma.

\begin{lemma}\label{proj}
Let $s>3$ and let $P_{i,j}$ be a sequence of probabilities, such that
\begin{equation*}
\limsup_{i,j\to\infty}{\frac{P_{i,j}}{(i+j)^{-s}}}<\infty.
\end{equation*}
Consider a shift invariant percolation model on $\Z^2$ on which a bond between
$(x_1,y_1)$ and $(x_2,y_2)$ is open with marginal probability
$P_{|x_1-x_2|,|y_1-y_2|}$. Assign conductance $1$ to every open bond,
and $0$ to every closed one. Call this electrical network $G_1$. Now,
perform the following projection process: for every open long
(i.e. not nearest neighbor)
bond $(x_1,y_1),(x_2,y_2)$ we
\\(A) erase the bond, and
\\(B) to each nearest neighbor bond in
$[(x_1,y_1),(x_1,y_2)]\cup[(x_1,y_2),(x_2,y_2)]$ increase the conductance by
$|x_1-x_2|+|y_1-y_2|$.
\\We call this new electrical network $G_2$. Then
\\(I) A.s. all of the conductances in $G_2$ are finite.
\\(II) The effective conductance of $G_2$ is bigger or equal to that of $G_1$.
\\(III) The distribution of the conductance of an edge in $G_2$ is shift
invariant.
\\(IV) If $s>4$ then the conductance of an edge is in $L^1$.
\\(V) If $s=4$ then the conductance $C_e$ of an edge has a {\em Cauchy tail},
i.e. there is a constant $\chi$ such that $\prob(C_e>n\chi)\leq n^{-1}$
for every $n$.
\end{lemma}

To complete the picture, we need the following theorem about random  electrical
networks on $\Z^2$. The theorem is proved in the next section.

\begin{theorem}\label{cauchytail}
Let $G$ be a random electrical network on the
nearest neighbor bonds of the lattice $\Z^2$, such that all of
the edges have the same conductance distribution, and this distribution has a
Cauchy tail. Then, a.s., a random walk on $G$ is recurrent.
\end{theorem}

Lemma \ref{proj} and Theorem \ref{cauchytail} imply the following theorem:

\begin{theorem}\label{two_d_recur}
Let $s\geq 4$ and let $P_{i,j}$ be probabilities, such that
\begin{equation}\label{sgadol4}
\limsup_{i,j\to\infty}{\frac{P_{i,j}}{(i+j)^{-s}}}<\infty.
\end{equation}
Consider a shift invariant percolation model on $\Z^2$ on which the bond between
$(x_1,y_1)$ and $(x_2,y_2)$ is open with marginal probability
$P_{|x_1-x_2|,|y_1-y_2|}$.
If there exists an infinite cluster, then the random walk on this cluster is
recurrent.
\end{theorem}
\begin{proof}
The case $s=4$ follows directly from \ref{proj} and Theorem \ref{cauchytail}. For
the case $s>4$, notice that if (\ref{sgadol4}) holds for some $s>4$, then it holds
for $s=4$ too.
\end{proof}

\begin{proof}[proof of Lemma \ref{proj}]

(I): We calculate the expected number of bonds that are projected on the edge
$(x,y),(x,y+1)$:
W.l.o.g, the projected bond starts at some $(x,y_1\leq y)$, continues through
$(x,y_2\geq y+1)$, and ends at some $(x_1,y_2)$. The expected number will be
\begin{eqnarray*}
2\sum_{y_1\leq y,y_2\geq y+1,x_1}{P_{|y_2-y_1|,|x_1-x|}}
&\leq& 4M\sum_{j\leq 0,k\geq 1,h\geq 0}{(k-j+h)^{-s}}\\
&\leq& 4M\sum_{l>0,h\geq 0}{(l+h)^{1-s}}\\
&\leq& 4M\sum_{l>0}{(l)^{2-s}}<\infty,
\end{eqnarray*}
where
\begin{equation*}
M=\sup_{i,j}{\frac{P_{i,j}}{(i+j)^{-s}}}<\infty.
\end{equation*}
and therefore (I) is true.
\\(II) let $E$ be a bond which is projected on a path of length $n$. $E$ has
conductance $1$, and is therefore equivalent to a sequence of $n$ edges with
conductance $n$ each. So, Divide $E$ that way. By identifying the
endpoints of these edges with actual vertices of the lattice, we only increase
the effective conductance of the network.
\\(III) is trivial.
\\(IV) and (V) follow from the same calculation performed in the proof of (I).
\end{proof}

\section{Random electrical networks}\label{elnet}

In this section we discuss random electrical networks.
We have two main goals in
this section:
\\{\bf Theorem \ref{cauchytail}.}
{\em Let $G$ be a random electrical
network on the lattice $\Z^2$, such that all
of the edges have the same conductance distribution, and this distribution has
a Cauchy tail. (Notice that we do not require any independence).
Then, a.s., a random walk on $G$ is recurrent.  }

and
\\{\bf Theorem \ref{iid}.}
{\em Let $G$ be a recurrent graph with bounded degree. Assign i.i.d.
conductances on the edges of $G$. Then, a.s., the resulting electrical
network is
recurrent.  }

Notice that if in Theorem \ref{cauchytail} we don't require a Cauchy tail,
then the network might be transient. A good example would be the projected
two-dimensional long-range percolation with $3<s<4$ (See Lemma \ref{proj}).

{}First, we prove Theorem \ref{cauchytail}, which is important for the
previous section. We need the following lemma, which sets some bound for the
sum of random variables with a Cauchy tail:

\begin{lemma}\label{ergod}
Let $\{f_i\}_{i=1}^\infty$ be identically distributed positive random variables
that have a Cauchy tail. Then, every $\epsilon$ has $K$ and $N$ such that
if $n>N$, then
\begin{equation*}
\prob\left(\frac{1}{n}\sum_{i=0}^{n}{f_i}>K\log n\right)<\epsilon.
\end{equation*}
\end{lemma}
\begin{proof}
$f_i$ has a Cauchy tail, so there exists $C$ such that for every $n$,
\begin{equation*}
\prob(f_i>n)<\frac{C}{n}.
\end{equation*}
Let $M>\frac{2}{\epsilon}$ be a large
number. Let $N$ be large enough that $CN^{1-M}<\halb\epsilon$. Choose $n>N$, and
let $g_i=\min (f_i,n^M)$ for all $1\leq i\leq n$. Then,
\begin{eqnarray*}
\prob\left(\frac{1}{n}\sum_{i=1}^{n}{f_i}
\neq\frac{1}{n}\sum_{i=1}^{n}{g_i}\right)
&\leq& n\cdot\prob(f_1\neq g_1)\\
\leq Cn^{1-M}&<&\halb\epsilon.
\end{eqnarray*}
$\E (g_i)\leq CM\log n$, and $g_i$ is positive. Therefore, by Markov's inequality,
if we take $K=CM^2$, then
\begin{equation*}
\prob\left(\frac{1}{n}\sum_{i=1}^{n}{g_i}>K\log n\right)<
\frac{CM\log n}{CM^2\log n}=\frac{1}{M}<\halb\epsilon.
\end{equation*}
and so
\begin{equation*}
\prob\left(\frac{1}{n}\sum_{i=1}^{n}{f_i}>K\log n\right)<
\epsilon.
\end{equation*}
\end{proof}
We use another lemma:
\begin{lemma}\label{1on}
Let $A_n$ be a sequence of events such that $\prob(A_n)>1-\epsilon$ for every
$n$, and let $\{a_n\}_{n=1}^\infty$ be a sequence s.t.
\begin{equation*}
\sum_{n=1}^\infty{a_n}=\infty.
\end{equation*}
Then, with probability at least $1-\epsilon$,
\begin{equation*}
\sum_{n=1}^{\infty}{1_{A_n}}\cdot{a_n}=\infty.
\end{equation*}
\end{lemma}
\begin{proof}
It is enough to show that for any $M$,
\begin{equation*}
\prob\left(\sum_{n=1}^{\infty}{1_{A_n}}\cdot{a_n}<M\right)\leq\epsilon.
\end{equation*}
Assume that for some $M$ this is false. Define $B_M$ to be the event
\begin{equation*}
B_M=\left\{\sum_{n=1}^{\infty}{1_{A_n}}\cdot{a_n}<M\right\}.
\end{equation*}
Since $\prob(B_M)>\epsilon$, we know that there exists $\delta>0$ such that
$\prob(A_n|B_M)>\delta$ for all $n$.
Therefore,
\begin{equation*}
\E \left(\sum_{n=1}^{\infty}{1_{A_n}}\cdot{a_n}|B_M\right)
\geq\delta\sum_{n=1}^{\infty}{a_n}=\infty,
\end{equation*}
which contradicts the definition of $B_M$.
\end{proof}
Now, we can prove Theorem \ref{cauchytail}.

\begin{proof}[Proof of theorem \ref{cauchytail}]
Let $G$ be a random electrical network on the lattice $\Z^2$, such that all of
the edges have the same conductance distribution, and this distribution has a
Cauchy tail.
\\Define the cutset $\Pi_n$ to be the set of edges exiting the
square $[-n,n]\cross[-n,n]$. We want to estimate
\begin{equation*}
\sum_n{C_{\Pi_n}^{-1}}.
\end{equation*}
Let $\epsilon>0$ be arbitrary. Let $e_n(i)$ be the $i$-th edge (out of
$(8n+4)$) in $\Pi_n$. By Lemma \ref{ergod}, there exist $K$ and $N$,
such that for every $n>N$, we have
\begin{equation}\label{casum}
\prob\left(\sum_{i=1}^{8n+4}{C(e_n(i))}\leq Kn\log n\right)>1-\epsilon.
\end{equation}
Call the event in equation (\ref{casum}) $A_n$.
Set ${a_n=(Kn\log n)^{-1}}$ for $n=N,...,\infty$.
Now,
\begin{equation*}
\sum_n{C_{\Pi_n}^{-1}}\geq\sum_{n=N}^\infty{1_{A_n}\cdot{a_n}}.
\end{equation*}
By the definition of $\{a_n\}$,
\begin{equation*}
\sum_{n=N}^{\infty}a_n=\infty.
\end{equation*}
On the other hand, $\prob(A_n)>1-\epsilon$ for all $n$.
So, by Lemma \ref{1on},
\begin{equation*}
\prob\left(\sum_n{C_{\Pi_n}^{-1}}=\infty\right)\geq 1-\epsilon.
\end{equation*}
Since $\epsilon$ is arbitrary, we get that a.s.
\begin{equation*}
\sum_n{C_{\Pi_n}^{-1}}=\infty.
\end{equation*}
\end{proof}

Now, we turn to prove Theorem \ref{iid}. First, we need a lemma:
\begin{lemma}[Yuval Peres]\label{yuval}
Let $G$ be a recurrent graph, and let $C_e$ be random conductances on the edges
of $G$. Suppose that there exists $M$ such that $\E (C_e)<M$ for each edge $e$.
Then, a.s., $G$ with the conductances $\{C_e\}$ is a recurrent electrical
network.
\end{lemma}
\begin{proof}
Let $v_0\in G$, and let $\{G_n\}$ be an increasing sequence of finite sub-graphs of $G$,
s.t. $v_0\in G_n$ for every $n$ and s.t.
$G=\cup_{n=1}^{\infty}G_n$. By the definition of effective conductance,
\begin{equation*}
\lim_{n\to\infty}C_{\ef}(G_n)=C_{\ef}(G)=0.
\end{equation*}
Let $X_n$ be the space of functions $f$ s.t. $f(v_0)=1$ and
$f(u)=0$ for every $u\in G-G_n$. We know that
\begin{equation*}
C_{\ef}(G_n)=\inf_{f\in X_n}\sum_{(v,w) \text{ is an edge in } G}{(f(v)-f(w))^2}.
\end{equation*}
If we denote by $H$ (resp. $H_n$) the electrical network of the graph $G$
(resp. $G_n$)
and conductances $C_e$, then
\begin{equation*}
C_{\ef}(H_n)=\inf_{f\in X_n}\sum_{(v,w) \text{ is an edge in G}}
{C_{v,w}(f(v)-f(w))^2}
\end{equation*}

Let $f\in X_n$. Denote $G_n(f)$ for
\begin{equation*}
\sum_{(v,w) \text{ is an edge in $G$}}{(f(v)-f(w))^2}
\end{equation*}
and $H_n(f)$ for
\begin{equation*}
\sum_{(v, w)\text{ is an edge in $G$}}{C_{v,w}(f(v)-f(w))^2}.
\end{equation*}
There exists an $f\in X_n$ such that $G_n(f)=C_{\ef}(G_n)$. Since
$\E (H_n(f))<MG_n(f)$,
we get that
\begin{equation*}
\E (C_{\ef}(H_n))\leq M(C_{\ef}(G_n).
\end{equation*}
So, by Fatou's lemma,
\begin{equation*}
\E(C_{\ef}(H))\leq\lim_{n\to\infty}\E(C_{\ef}(H_n))\leq
M\lim_{n\to\infty}(C_{\ef}(G_n))=0,
\end{equation*}
and therefore $C_{\ef}(H)=0$ a.s.
\end{proof}

Now we can prove Theorem \ref{iid}. The main idea is to change the conductances
in a manner that will not decrease the effective conductance, but after this
change, the conductances will have bounded expectations (although they might
be dependent).
\begin{proof}[Proof to Theorem \ref{iid}]
Let $G$ be a recurrent graph, and let $d$ be the maximal degree in $G$.
Let $\{C_e\}_{\{e \text{ is an edge in }G\}}$ be i.i.d.
non-negative variables, and let $H$ be the electrical network defined on the
graph $G$ with the conductances $\{C_e\}$. We want to prove that with
probability one $H$ is recurrent. Let $M$ be so large that
\begin{equation*}
P(C_e\geq M)<\frac{1}{d^5}.
\end{equation*}
We introduce some notation: edges whose conductances are bigger than $M$
will be called $bad$ edges. Vertices which belong to bad edges will also be
called bad. We look at connected clusters of bad edges. Edges that are good
but have at least one bad vertex, will be called {\em boundary} edges.

By the choice of $M$, the sizes of the clusters of bad edges are dominated by
sub-critical Galton-Watson trees. Define a new network $H'$ as
follows: Let $U(e)$ be the connected component to which $e$
belongs (if $e$ is bad) or to which $e$ is attached (if it is a boundary edge).
If $e$ is in the boundary of two components, then we take $U(e)$ to be their
union. For a bad or boundary $e$, the new conductance will be
$2M\cdot (\#U(e)+\#\partial U(e))^2$, where $\#$ measures the number of edges.
If $e$ is a good edge then its conductance is unchanged.
The size of the connected cluster satisfies
\begin{equation*}
\prob(\#U(e)+\#\partial U(e)>n)=o(n^{-4}).
\end{equation*}
Therefore, the expected values of the conductances of the edges are
uniformly bounded,
So by Lemma
\ref{yuval}
$H'$ is recurrent. All we need to prove is that the effective resistance of
$H'$ is not bigger than that of $H$: Let $F$ be a flow, and let $U$ be a
connected component of bad edges in $G$.
The energy of $F$ on $U$ in the network $H$ will
be
\begin{equation*}
\en_{U,F}(H)=\sum_{e\in U\cup \partial U}{\frac{F_e^2}{C_e}}\geq
\sum_{e\in \partial U}{\frac{F_e^2}{C_e}}\geq
\sum_{e\in \partial U}{\frac{F_e^2}{M}}.
\end{equation*}
{}For every $e$ in $U\cup \partial U$, the flow $F_e$ is smaller than
\begin{equation*}
\sum_{e'\in \partial U}|F_{e'}|,
\end{equation*}
so
\begin{equation*}
F^2_e\leq \#\partial U\cdot\sum_{e'\in \partial U}{F_{e'}^2}
\leq M\cdot\#\partial U\cdot \en_{U,F}(H).
\end{equation*}
Therefore,
\begin{eqnarray*}
\en_{U,F}(H')&=&
\sum_{e\in U\cup \partial U}{\frac{F_e^2}{2M\cdot (\#U+\#\partial U)^2}}\\
&\leq&(\#U+\#\partial U)\frac{M\cdot\#\partial U\cdot \en_{U,F}(H)}
{2M\cdot (\#U+\#\partial U)^2}\leq \en_{U,F}(H).
\end{eqnarray*}
Thus, by Thomson's theorem (see \cite{yuval}), the effective resistance of
$H'$ is not bigger than that of $H$, and we are done.
\end{proof}

\section{Critical behavior of the free long-range
random cluster model}\label{frcm}
We return to the critical behavior. Our goal in this section is to prove Theorem
\ref{fr_int} and Corollary \ref{ising_extr}
We begin with the following extension of Theorem \ref{fr_int}:
\begin{theorem}\label{fcr}
Let $d<s<2d$ and let $\{P_k\}_{k\in\Z^d}$ be nonnegative numbers such that
$\forall_k(P_k=P_{-k})$ and
\begin{equation}\label{smas}
\liminf_{\|k\|\to\infty}\frac{P_k}{{\|k\|_1^{-s}}}>0.
\end{equation}
Let $\beta>0$, and consider the infinite volume limit of the free random
cluster model with probabilities $1-e^{-\beta P_k}$ and with $q\geq 1$ states.
Then, a.s., at
\begin{equation*}
\beta_c=\inf(\beta |\text{ a.s. there exists an infinite cluster})
\end{equation*}
there is no infinite cluster.
\end{theorem}
We need the following extension of Lemma \ref{normal1}:
\begin{lemma}\label{normal2}
Let $d\geq 1$.
Consider an ergodic (not necessarily independent) percolation model on $\Z^d$
which satisfies
\begin{equation}\label{dominate}
\prob\left(i\leftrightarrow j | \borel_{i,j}\right)\geq P_{i-j}
\end{equation}
Where $i\leftrightarrow j$ denotes the event of having an open bond
between $i$ and $j$, and $\borel_{i,j}$ is the $\sigma$-field created by all
of the events $\{i'\leftrightarrow j'\}_{(i',j')\neq(i,j)}$.
Assume further that:
\\(A) The distribution has the FKG property {\rm \cite{FKG}}.
\\(B) A.s. there is a unique infinite cluster.
\\(C) There exists $d<s<2d$ s.t.
\begin{equation*}
\liminf_{||k||\to\infty}\frac{P_k}{||k||^{-s}}>0.
\end{equation*}
Then, for every $\epsilon>0$ and $\rho$ there exists $N$ such that with
probability bigger than $1-\epsilon$, inside the cube $[0,N-1]^d$ there exists
an open cluster which contains at least $\rho N^{\frac{s}{2}}$ vertices.
\end{lemma}
Lemma \ref{normal2} is proved exactly the same way as Lemma \ref{normal1}. Lemma
\ref{normal2} is valid for the free random cluster model measure considered in
Theorem \ref{fcr}. We can use Lemma \ref{normal2} to prove the
following:
\begin{lemma}\label{normalf}
Let $d<s<2d$ and let $\{P_k\}_{k\in\Z^d}$ be nonnegative numbers such that
$\forall_k(P_k=P_{-k})$ and
\begin{equation}\label{smals}
\liminf_{|k|\to\infty}\frac{P_k}{{\|k\|_1^{-s}}}>0.
\end{equation}
Let $\beta>0$, and consider the infinite volume limit of the free random
cluster model with probabilities $1-e^{-\beta P_k}$. Assume that, a.s., there
is an infinite cluster.
Then, for every $\epsilon$ and $\rho$ there is an $N$ such that given the
the values (open or closed) of all of the edges that have at least one end
point out of the cube $[0,N-1]^d$, the probability of having an open cluster
of size $\rho N^{\halb s}$ within $[0,N-1]^d$ is larger than $1-\epsilon$.
\end{lemma}
\begin{proof}
The proof follows the guideline of the proof of Lemma \ref{normal1}: Choose
$\epsilon'$ and $\theta$, and let $M$ be s.t. by Lemma \ref{normal2} with
probability larger than $1-\epsilon'$ there exists an open cluster of size
$\sqrt{\theta}M^{\halb s}$ inside $[0,M-1]^d$. Let $K$ be s.t. this probability is larger than
$1-2\epsilon'$ even if all of the edges with at least one endpoint out of
$[-K,K+M-1]^d$ are closed. Such $K$ exists because the free measure on $\Z^d$ is
the limit of the free measures on $[-K,K+M-1]^d$
when $K$ tends to infinity. Now, let $R$ be a large number.
Assume that all of the edges with (at least) one
endpoint out of $[-K,RM+K-1]^d$ are closed. For a cube
\begin{equation*}
\C=\prod_{j=1}^d{[l_iM,(l_i+1)M-1]}\qquad\qquad 0\leq l_i\leq R-1
\end{equation*}
in $[-K,RM+K-1]^d$, the probability of the cube to be {\em alive}, i.e. to have
an open cluster of size $\sqrt{\theta}M^{\halb s}$
is larger than $1-2\epsilon'$
(because of domination). The probability that there exists an open bond between two
living cubes that are $k$ cubes away from each other is larger than
$\eta_s(\frac{\theta}{2},k)$.
Now, we can proceed exactly as in the proof of Lemma
\ref{normal1}. With $\epsilon'$, $\theta$ and $R$ properly chosen, the lemma
is proved.
\end{proof}
Now, we can prove Theorem \ref{fcr}:
\begin{proof}[Proof of Theorem \ref{fcr}]
Let $\{P_k\}_{k\in\Z^d}$ be such that for every $k$, $P_k=P_{-k}$ and such that
\begin{equation*}
\kappa=\liminf_{\|k\|\to\infinity}\frac{P_k}{\|k\|_1^{-s}}>0.
\end{equation*}
Let $\beta$ be s.t. for the Random Cluster Model with interactions $\{P_k\}$
and inverse temperature $\beta$ there exists, a.s., an
infinite cluster. What we need to show is that there exists an $\epsilon>0$
s.t. there exists an infinite cluster at inverse temperature $\beta-\epsilon$.
For every $a$ and $b$ consider the independent percolation model $\idp(a,b,s)$
where every vertex exists with probability $a$ and two vertices $x$ and $y$
are attached to each other with probability $1-e^{-b|x-y|^{-s}}$.
Let $\gamma$, $\lambda$ and $\delta$ be s.t. in
$\idp(\lambda-\delta,\gamma-\delta,s)$ there exists, a.s., an infinite cluster.

Let $N$ be so large that by Lemma \ref{normalf} with probability larger than
$\lambda$ there exists a cluster of size $\rho N^{\halb s}$ inside
$[0,N-1]^d$, where the probability is with respect to the free measure on
$[0,N-1]^d$, and $\rho$ is s.t.
\begin{equation}\label{rho}
\frac{\rho^2}{2q}>\gamma.
\end{equation}
By the choice of $\rho$ (\ref{rho}) we get that
the probability of having an open bond between clusters of size
$\rho N^{\halb s}$ that are located in the cubes at $Nx$ and $Ny$ is
(no matter what happens in any other bond)
at least $1-e^{-\gamma\|x-y\|_1^{-s}}$.

Now, let $\epsilon>0$ be s.t. in inverse temperature $\beta-\epsilon$ the
probability of having this big cluster is larger than $\lambda-\delta$, and the
probability of having an open bond is larger than $e^{(\gamma-\delta)|x-y|^{-s}}$.
Such $\epsilon$ exists, because the probability of any event in a finite
random cluster model is a continuous function of the (inverse) temperature.
When considering the renormalized model in inverse temperature
$\beta-\epsilon$, it dominates
$\idp(\lambda-\delta,\gamma-\delta,s)$, and therefore has an infinite cluster.
\end{proof}
We can now restate and prove Corollary \ref{ising_extr}:
\\{\bf Corollary \ref{ising_extr}.}
{\em
Let $\{P_k\}_{k\in\Z^d}$ be nonnegative numbers s.t. $P_k=P_{-k}$ for every $k$
and s.t. $P_k\sim \|k\|_1^{-s}$ ($d<s<2d$).
Consider the Potts model with $q$ states on $\Z^d$, s.t. the interaction between
$v$ and $u$ is  $P_{v-u}$. At the critical temperature,
the free measure is extremal.
}
\begin{proof}[Proof of Corollary \ref{ising_extr}]
Recall the following construction of a configuration of the free measure
of the Potts model:
choose a configuration of the free measure of the Random Cluster model, and
color each of the clusters by one of the $q$ states. The states of different clusters are
independent of each other.
By Theorem \ref{fr_int}, there is no infinite cluster at the critical temperature.
Therefore, for every $n$ and $\epsilon$ there exists $K$ s.t. with probability $1-\epsilon$
for every $x$ s.t. $\|x\|_1\leq n$ and $y$ s.t. $\|y\|_1\geq K$, $x$ and $y$ belong to distinct
clusters.

Therefore, for the Potts model, there is an event $E$ of probability bigger than
$1-\epsilon$ s.t. given $E$, the coloring of $\{x:\|x\|_1\leq n\}$ is independent of the
coloring of $\{y:\|y\|_1\geq K\}$. Therefore, the tail $\sigma$-field
\begin{equation*}
\bigcap_{K=1}^{\infty}\sigma\left(v \text{ s.t. } \|v\|_1>K \right)
\end{equation*}
is trivial, and therefore the measure is extremal.
\end{proof}
\section{Remarks and problems}
Many more questions can be asked about these clusters. One example is the
volume growth rate. It can be shown that the growth of the infinite cluster
is not bigger than exponential with the constant
\begin{equation*}
\sum_{k\in\Z^d}P_k.
\end{equation*}

In the case $d<s<2d$, The growth can be bounded from below by $\exp(n^{\phi(s)})$, for
$\phi(s)=\log_2(2d/s)-\epsilon$.
This can be proved as follows: if $\beta$ is large enough, then in the proof of
Theorem \ref{thm:ns}, we may take $C_n=\exp(2^{\phi(s)\cdot n})$. Then, the
$n$-th degree cluster contains
\begin{equation*}
\prod_{k<n}{C_k},
\end{equation*}
vertices, while its diameter is at most $2^n$. This gives a lower bound of
$\exp(n^{\phi(s)})$ for the growth. if $\beta$ is not so large, then by using
Lemma \ref{normal1} we can make it large enough.
\\In the case $s=2d$, the volume growth rate is subexponential
(see \cite{BenBer}). In the case $s<2d$ it is not known. So, we get
a few questions on the structure of the infinite cluster.

\begin{question}\label{growth}
What is the volume growth rate of the infinite clusters of super-critical
long-range percolation with $d<s<2d$?
Is it exponential?
\end{question}
\begin{question}
How many times do two independent random walks paths on the infinite cluster of
long-range percolation intersect?
\end{question}
\begin{question}
Are there any nontrivial harmonic functions on the infinite cluster of
one-dimensional long-range percolation with $d<s<2d$?
\end{question}
Other questions can be asked on the critical behavior. The renormalization
lemma (Lemma \ref{normal1}) is only valid when $d<s<2d$. So, the arguments
given here say nothing about the critical behavior on other cases. At the case
$d=1$ and $s=2$, Aizenman and Newman proved that there exists an infinite
cluster at criticality (see \cite{AN}). For the other cases the following
questions are still open:
\begin{question}\label{sgeq2d}
Does critical long-range percolation have an infinite cluster when
$d \geq 2$ and
$s\geq 2d$?
\end{question}
As remarked by G. Slade, the methods used in \cite{hasl} might be used to prove
that for $d>6$ and $s>d+2$ there is no infinite at criticality.
This can reduce Question \ref{sgeq2d}
to the case $2\leq d\leq 6$.

\begin{question}
Does the conclusion of Theorem \ref{jeff} hold for sequences which decay
faster than those treated in Theorem \ref{jeff} and slower than those treated
by Steif and Meester
{\it {\em (\cite{steif})}}? i.e.

Let $d\geq 2$. For which percolating $d$-dimensional arrays of probabilities
$\{P_k\}_{k\in \Z^d}$ there exist an $N$ s.t. the independent
percolation model with probabilities
\begin{equation*}
P'_k = \left\{
        \begin{array}{ll}
                P_k             &\|k\|_1<N\\
                0               &\|k\|_1\geq N
        \end{array}
        \right.
\end{equation*}
also has, a.s., an infinite cluster?
\end{question}
The arguments given in this paper are not strong enough to prove that there is
no infinite cluster in the wired random cluster model at the critical
temperature. So, the following question is still open:
\begin{question}\label{wired}
Is there an infinite cluster at the critical temperature in the wired random
cluster model with $d<s<2d$?
\end{question}
A different formulation of the same question is
\\{\bf Question \ref{wired} (Revised).}
{\em Let $d\geq 1$ and let $d<s<2d$. Let $\{P_k\}_{k\in Z^d}$ be s.t. $P_k=P_{-k}$
for every $k$ and s.t. $P_k\sim \|k\|_1^{-s}$. Consider the Potts model (with $q$ states)
with interaction
$P_{u-v}$ between $u$ and $v$. Let $\beta$ be the critical inverse temperature for
this Ising model. Is there a unique Gibbs measure at inverse temperature $\beta$?
}

Question \ref{wired} is related to the question whether the free and the wired
measures agree on the critical point. Conjecturing that for high values of $q$,
the number of states,
the critical wired measure has an infinite cluster, we will get the conjecture
that the two measures won't agree at the critical point.


\appendix
\section{Proof of Lemma \ref{lem:aldous}}\label{app:aldous}
The proof of Lemma \ref{lem:aldous} is based on the methods from \cite{aldous}. As in \cite{aldous}, we construct simultaneously
a random walk and the random graph. We cite a result from \cite{aldous} on the connection between excursions of the random
walk and the connected components of the graph. We then use this result to prove Lemma \ref{lem:aldous}. 

The construction of the random walk and the random graph is as follows: For each ordered pair $(i,j),i\neq j$, let $U_{i,j}$ be an exponential 
$M^{-\xi}m(h(j))$ variable, independent over pairs. Choose $v_1$ by size biased sampling (i.e. the probability that $v_1=h_i$ is proportional to
$m(h_i)$). Let $\{v:U_{v_1,v}\leq m(v_1\}$ be the set of children of $v_1$, and order them as $v(2),v(3),\ldots$ so that $U_{v_1,v(i)}$ is increasing.
Start the walk $z(\cdot)$ with $z(0)=0$, and let 
\[
z(u)=-u+\sum_v m(v){\bf 1}_{(U_{v_1,v}\leq u)}, \ \ \ \ \ \ \ \ 0\leq u\leq m(v_1).
\]
In particular,
\[
z(m(v_1)) = -m(v_1) +\sum_{v\mbox{ child of } v_1}m(v).
\]
Inductively, write $\tau_{i-1}=\sum_{j\leq i-1}m(v_j)$. If $v_i$ is in the same component as $v_1$, then the set 
\[
\{v\notin\{v_1,\ldots,v_{i-1}\}: v\mbox{ is a child of one of } \{v_1,\ldots,v_{i-1}\} \}
\]
consists of $v_1,\ldots,v_{l(i)}$ for some $\l(i)\geq i$. Let the children of $v_i$ be 
\[
\{v\notin\{v_1,\ldots,v_{l(i)}\}: U_{v_i,v} \leq m(v_i)\},
\]
and order them as $v_{l(i)+1},v_{l(i)+2},\ldots$ such that $U_{v_i,v}$ is increasing.
Set
\begin{equation*}
z(\tau_{i-1}+u)=z(\tau_{i-1})-u + \sum_{v \mbox{ child of } v_i} m(v){\bf 1}_{U_{v_i,v}<u}, \ \ \ \ \ \ \ 0\leq u \leq m(v_i).
\end{equation*}
After exhausting the component containing $v_1$, choose the next vertex by size biased sampling among the remaining vertices. Continue.
For simplicity, for $u>M$ we define $z(u)=z(M)+M-u.$

This construction yields a forest on the vertices $h_1,\ldots,h_k$, an ordering $v_1,\ldots,v_k$ and a walk $z(u); \ 0\leq u \leq M$. Add extra edges
between $v_i$ and $v_j$ for every pair such that $i<j\leq l(j)$ and $U_{v_i,v_j}\leq m(v_i)$. The resulting random graph has the same distribution as the
inhomogenous random graph, the ordering of the vertices $v_1,\ldots,v_k$ is size biased, and the relation between the components of the graph and the
random walk is as appears in the lemma below:
\begin{lemma}[\cite{aldous}, Page 828]\label{lem:appfromaldous}
Every connected component in the graph is a sequence of vertices ${v_i,v_{i+1}\ldots,v_j}$ such that 
\[
z(\tau_j)=z(\tau_i)-m(v_i), \ \ \ z(u)\geq z(\tau_j) \mbox{ on } \tau_{i-1} < u < \tau_j.
\]
Furthermore, the size of the component is $\tau_j-\tau_{i-1}$.
\end{lemma}

We now use this construction to prove Lemma \ref{lem:aldous}.

For $0\leq u\leq M$, define $i(u)=\min\{i:\tau_{i}\geq u\}$ to be the particle that is being processed at time $u$.
We define the set $B(u)$ to be the set of all particles seen up to time $u$, namely
\[
B(u)=\{v_j:j\leq i(u) \mbox{ or } \exists_{k<i}\mbox{ s.t. $j$ is a child of $k$ } \mbox{ or } U_{i(u),j}<u-\tau_{i(u)-1}\}.
\]
Then define the {\em drift} $D(u)$ to be
\[
D(u)= - 1 + M^{-\xi}\sum_{h_i\notin B(u)}m^2(h_i).
\]
$D(\cdot)$ is the drift of $z(\cdot)$ in the sense that
\[
I(u):=z(u)-\int_0^u D(s)ds
\]
is a martingale.

Clearly, $D(u)$ is decreasing with $u$.

We remember that $1<\xi<2$, and that $\gamma$ is chosen such that $1>\gamma>\frac{16+\xi}{18}$. We also take $\gamma'$ s.t. $\gamma'>\frac{4+\xi}{6}$ and $3\gamma - 2 > \gamma'$.
\ignore{
We take $\epsilon$ so that 
$\epsilon<\gamma-\gamma'$ and $\epsilon<(2\gamma'-\xi)/4$, but $\epsilon>1-\gamma$.
}
Note that $\gamma-\gamma'>2(1-\gamma)$, and take $\gamma-\gamma'> \epsilon >2(1-\gamma)$.
Then, $\epsilon/2>1-\gamma$, $\epsilon<\gamma-\gamma'$ and $\epsilon<(2\gamma'-\xi)/4$.
Let $\alpha=\xi-\gamma'-\epsilon$. Note that $\alpha-\gamma'=\xi-2\gamma'-\epsilon<-5\epsilon$. Let $\delta<\epsilon/2-(1-\gamma)$, and let $\theta$ be so that $\gamma'-\delta<\theta<\gamma'$.

\ignore{
\[
\gamma-\gamma' < (2\gamma'-\xi)/4 ?
\]

\begin{eqnarray*}
\gamma-\gamma' < 1-\gamma'  < 1 - \frac{4+\xi}{6} = \frac{2-\xi}{6}< (2\gamma'-\xi)/4
\end{eqnarray*}
}

\begin{claim}\label{claim:bigpart}
For $M$ large enough, with probability larger than $1-e^{-M^{(\gamma-\gamma'-\epsilon)/2}}$, for every $u>\frac{1}{2}M^\gamma$, and every $h_i\notin B(u)$, we have
$m(h_i)<M^\alpha$.
\end{claim}

\begin{proof}
As $B$ is increasing, it is enough to speak about $u=\frac{1}{2}M^\gamma$. By the construction, for a given particle $h_i$ with $m(h_i)\geq M^\alpha$,
\[
P(h_i\notin B(u))\leq \exp(-u M^{-\xi} M^\alpha)
=\exp\left(-\frac 12 M^{\gamma + \alpha -\xi}\right) = \exp\left(-\frac 12 M^{\gamma-\gamma'-\epsilon}\right),
\]
and since $\gamma-\gamma'-\epsilon>0$ and there are at most $M^{1-\alpha}<M$ such particles, the claim follows.
\end{proof}


Let $A$ be the event $A=\{\forall_{h_i\notin B(M^\gamma/2)}\ m(h_i)<M^\alpha\}$. Then by the previous claim $P(A)\geq 1-e^{-M^{(\gamma-\gamma'-\epsilon)/2}}$.

We calculate the variance of the difference of $I(\cdot)$ within one time unit. Assume that $u>\frac 12M^\gamma$. Let $\F_u$ be the $\sigma$-algebra generated by the process up to time $u$.
Remembering that 
$E(I(u+1)|\F_u)=I(u)$, we get
\begin{eqnarray}\label{eq:varmart}
\nonumber
\var (I(u+1)-I(u)|\F_u ; A) &\leq& 
(D(u)+1)^2 + M^{-\xi}\sum_{i\notin B(u)} m^3(h_i)\\
&\leq& (D(u)+1)(D(u)+1+M^{\alpha}).
\end{eqnarray}

We now divert our attention to the rate at which $D(u)$ decreases after time $M^\gamma/2$. Let 
\[
L(u)=D(u)+1=M^{-\xi}\sum_{i\notin B(u)}m^2(h_i).
\]

\begin{lemma}\label{lem:decdrift} Let $\kappa=\frac 13 M^{\gamma'}$. Let $u> \frac 12M^\gamma$. Let $L:=L(u)$. Let 
\[
\ell=\left\lceil\log_2\left(2M^\alpha/LM^{\xi-1}\right)\right\rceil.
\]
Then conditioned on the event $A$,
\begin{equation}\label{eq:decdrift}
P\left(\left.L(u) - L(u+\kappa)  < \frac{L}{16\ell}\cdot \left[1-e^{-\frac L2 \frac\kappa M}\right] \right| \F_u ; A\right)
\leq \exp\left(-  \left[1-e^{-\frac L2 \frac\kappa M} \right] \cdot   \frac{LM^\xi}{32\ell M^{2\alpha}} \right).
\end{equation}
\end{lemma}
\begin{proof} 
For every $i\notin B(u)$, the (conditional) probability $P_i$ that $i\notin B(u+\kappa)$ satisfies
$
P_i\leq e^{-M^{-\xi}\kappa m(h_i)}.
$

Let $L=L(u)$.

Note that
\begin{eqnarray*}
M^{-\xi}\sum_{i:m(h_i) < \frac 12LM^{\xi-1}}m^2(h_i)
&\leq& \frac 12LM^{\xi-1} M^{-\xi}\sum_{i:m(h_i) < \frac 12LM^{\xi-1}}m(h_i)\\
&\leq& \frac 12LM^{\xi-1} M^{-\xi} M = \frac 12L.
\end{eqnarray*}

Therefore, 
\begin{eqnarray*}
M^{-\xi}\sum_{i:m(h_i) \geq \frac 12LM^{\xi-1}}m^2(h_i)
\geq \frac 12L.
\end{eqnarray*}

Recall that 
\[
\ell=\left\lceil\log_2\left(2M^\alpha/LM^{\xi-1}\right)\right\rceil.
\]
For $i=0,\ldots,\ell$, let 
\[
B_i=\left\{h_i\notin B(u):
2^i\leq \frac{m(h_i)}{\frac 12LM^{\xi-1}} <2^{i+1}
\right\}.
\]
Let $i_0$ maximize 
\[
\sum_{h_i\in B_i}m^2(h_i),
\] 
and let $B:=B_i$.
Then conditioned on $A$,
\begin{equation}\label{eq:inB}
M^{-\xi}\sum_{h_i\in B} m^2(h_i)\geq \frac{L}{2\ell}.
\end{equation}
$P(h_i\in B(u+\kappa))\geq 1-\exp\left(-\frac L2 \frac\kappa M\right)$ for each $h_i\in B$, and $|B|\geq \frac{LM^\xi}{2\ell M^{2\alpha}}$. Therefore by standard binomial estimates,
\begin{eqnarray*}
P\left(L(u) - L(u+\kappa)  < \frac{L}{16\ell}\cdot \left[1-e^{-\frac L2 \frac\kappa M}\right]  \right)
\leq \exp\left(-  \left[1-e^{-\frac L2 \frac\kappa M}\right] \cdot   \frac{LM^\xi}{32\ell M^{2\alpha}} \right)
\end{eqnarray*}

as desired.

\end{proof}

Note that if
$L$ is of constant order of magnitude or larger, then
the bound in \eqref{eq:decdrift} decays exponentially with a positive power of $M$. Thus, applying Lemma \ref{lem:decdrift} again and again, we get the following corollary regarding the decrease of the drift:


\begin{corollary}\label{cor:decdrift}
Let $L_0=L(M^\gamma/2)$.
Then with probability larger than $1-M^{-1}$, for every $k$ such that 
\[
L_0-\frac{kM^{\theta - 1}}{2}>\frac 12,
\]
we have
\begin{equation}\label{eq:decdrift2}
L(M^\gamma+k\kappa)\leq L_0-\frac {kM^{\theta-1}}{2}.
\end{equation}

Furthermore, with probability larger than $1-M^{-1}$, for every $u$ such that $L(u)>\frac 12$ and $u>\frac 12 M^\gamma$,
\begin{equation}\label{eq:decdrift2}
L(u+\kappa)\leq L(u)-\frac {M^{\theta-1}}{2} \ \ \ \mbox{ and } \ \ \ L(u-\kappa)\geq L(u)+\frac {M^{\theta-1}}{2}.
\end{equation}

\end{corollary}
Let $u_0=\max\left(\inf\{u:L(u)\leq 1
\},\ M^\gamma/2\right)$. Let $u_k=u_0+k\kappa$, $k\in\Z$. Let $T\leq\infty$ be the time at which the first excursion generated by a component larger than $M^\gamma$ ends.
If no such component exists, then $T=\infty$.

\begin{lemma}\label{lem:first_block}
There exists 
$\varphi>0$ such that for every $M$ large enough and every inhomogoneus random graph with size $M$ and parameter $\xi$, 
\[
P(T\geq u_0) > 1-M^{-\varphi}.
\]
\end{lemma}
\begin{proof}

Let $u\geq \frac 12M^\gamma$.
Calculation yields:
\[
E\left[\left.
\sum_{h_i\in B(u+\kappa)\setminus B(u)}m(h_i)
\right|\F_u\right]\geq (1-M^{-\epsilon})\kappa L(u).
\]
\ignore{
\[
P(h_i)=1-\exp(-M^{-\xi}m(h_i)\kappa)\geq M^{-\xi}m(h_i)\kappa(1 - M^{-\xi}m(h_i)\kappa)
\]
so we need
\[
M^{-\xi}m(h_i)\kappa<M^{-\epsilon}
\]
Check:
\[
M^{-\xi}m(h_i)\kappa < \frac 13 M^{-\xi+\alpha+\gamma'}<M^{-\epsilon}
\]

}
\ignore{
and
\[
\var\left[\left.
\sum_{h_i\in B(u+\kappa)\setminus B(u)}m(h_i)
\right|\F_u\right]\leq M^\alpha\kappa L(u).
\]

Let $X=\sum_{h_i\in B(u+\kappa)\setminus B(u)}m(h_i)$.
Then, if $L(u)\geq 1+3M^{-2\epsilon}$, by chebichef's inequality we get
\begin{eqnarray*}
P\left[\left.
z(u+\kappa)-z(u) \leq \kappa(L(u)-1)/2
\right|\F_u\right]
&=& 
P\left[\left.
X \leq \kappa\big[1+(D(u))/2\big]
\right|\F_u\right]\\
\leq
P\left[\left.
X \leq E(X|\F_u)  - \kappa D(u)/2
\right|\F_u\right]
&\leq&
\frac{\var (X|\F_u)}{\kappa^2 D^2(u)/4}
\leq 
\frac{1}{12} \frac{L(u) M^{\alpha - \gamma'}}{D^2(u)}.\\
&\leq&
\frac{1}{12} \frac{L(u) M^{-5\epsilon}}{D^2(u)}.
\end{eqnarray*}
}
Let $W(u)$ be the event $\{z(u+\kappa)-z(u) \geq \kappa D(u)/2\}$. Let $E$ be the event 
\[
E=
\left\{
\forall_{0\leq k \leq M^\gamma/2\kappa}
D(M^\gamma-k\kappa)\geq 
k\frac{M^{\theta-1}}{2}
\right\}.
\]
Note that by Corollary \ref{cor:decdrift}, $P(E^c\ ;\ T<u_0)<M^{-1}$. Therefore, it suffices to estimate $P(E\ ;\ T<u_0)$.
By the exponential Markov inequality, for every $u>M^\gamma/2$
\[
P\left[\left.
W^c(u);A;D(u)\geq 3M^{-\epsilon}
\right|\F_u\right] \leq \frac{\exp(-M^{2\epsilon})}{16}.
\]

\ignore{
If $D(u)> 3M^{-\epsilon}$, then $L(u)>1+3M^{-\epsilon}$ and $(1-M^{-\epsilon})L(u)>1+D(u)-1.1M^{-\epsilon}>1+D(u)/2$. In particular, the gap is of order of magnitude of $D(u)$. Now exponential Markov:
Checking exponential Markov for this case:

\begin{eqnarray*}
\sum P_iX_i = M^{\gamma'}\ ; \ \sup X_i = M^\alpha < M^{\gamma'-5\epsilon}\\
P\left(X< M^{\gamma'}(1-M^{-\epsilon})\right)\leq\frac{E(e^{-KX})}{e^{-KM^{\gamma'}(1-M^{-\epsilon})}}
\end{eqnarray*}

\[
\log\left(1-P_i+P_ie^{-KX_i}\right) < -P_i(KX_i - K^2X_i^2) < -P_iKX_i + M^\alpha K^2P_iX_i
\]
\[
\sum\log(...) < -KM^{\gamma'} + K^2M^\alpha M^{\gamma'}
\]

So we need to find $K$ that minimizes
\[
K^2M^{\alpha+\gamma'} - KM^{\gamma'-\epsilon}
\]
\begin{eqnarray*}
2KM^{\alpha+\gamma'} = M^{\gamma'-\epsilon}\\
K=M^{\gamma'-\epsilon-\gamma'-\alpha}/2=M^{-\epsilon-\alpha}/2
\end{eqnarray*}
\[
K^2M^{\alpha+\gamma'} - KM^{\gamma'-\epsilon} = M^{\gamma'-2\epsilon-\alpha}/4>M^{\epsilon}/4.
\]

need to check: for small $x$,
\[
\log(1-p+pe^{-x})<-px+px^2?
\]

For small $x$,
\[
\log(1-p+pe^{-x}) < \log(1-p+p(1-x+x^2)) = \log(1-px+px^2)
<p(-x+x^2)-p^2(cx^2)...
\]

So all of this works if $x$ is (objectively) small enough and $p\ll 1$, which is, indeed the case.
}

Let $k_0=3M^{1-\epsilon-\theta}$, and let $W$ be the event
\[
W=
\bigcap_{k>k_0:u_0-k\kappa\geq M^\gamma/2}
W(u_0-k\kappa).
\] 

Then $P(E\ ; W^c)<M^{-1}$. Therefore, it suffices to prove that if both the events $E$ and $W$ occur, then $T\geq u_0$. Under the event $E\cap W$, we have that 
\[
z(3M^\gamma/4) - z(M^\gamma/2) \geq
\frac{\kappa M^{\theta-1}}{4}\sum_{k=0}^{M^{\gamma}/4\kappa}\left(M^{\gamma}/2\kappa-k\right)
\geq \frac {M^{2\gamma+\theta-\gamma'-1}}{576}\geq M^{2\gamma-\delta-1}.
\]

Therefore, for under the event $E\cap W$, for every $u$ between $3M^\gamma/4$ and $u_0$ we have $z(u)>z(M^\gamma/2)$,
and therefore $u_0$ is in the same component excursion as $3M^\gamma/4$, and thus $T>u_0$.


\ignore{
between $u_0$ and $u_{-k_0}$ we have distance of $M^{1-\epsilon-\theta+\gamma'}$. So we need to show that
\[
1-\epsilon-\theta+\gamma' < 2\gamma-\delta-1.
\]
\begin{eqnarray*}
1-\epsilon-\theta+\gamma' < 1-\epsilon + \delta
\end{eqnarray*}
So we need 
\begin{eqnarray*}
2-\epsilon+2\delta < 2\gamma.
\end{eqnarray*}
i.e.
\begin{eqnarray*}
1-\gamma < \epsilon/2 - \delta
\end{eqnarray*}
}
\ignore{
To this end we start by
estimating the probability of the event $W(u)$ for various values of $u$.
Let $f(u)=M^{2\epsilon}D(u)$.

}

\end{proof}

\begin{lemma}\label{lem:other_blocks}
With probability at least $1-CM^{-\varphi}$,
there exists no excursion which is generated by a component larger than $M^\gamma$ which starts after time $u_0$.
\end{lemma}

\begin{proof}
The calculation is similar to the one from the previous proof. First, for $k=0,1,2,\ldots,\frac{M^{1-\theta}}{2}$,
\begin{eqnarray*}
&&P\left(
\exists_{u\in[u_{k+1},u_{k+2}]} I(u)>I(u_{k})+\kappa |D(u_{k})|/2
\right)\\
&\leq& \frac{4\kappa M^\alpha }{\kappa^2 D(u_{k})^2}
\leq 48M^{\alpha-\gamma'}\cdot M^{2(1-\theta)}k^{-2}\\
&=& Ck^{-2} M^{2+\xi-\epsilon-2\gamma'-2\theta} = C\frac{M^{-\varphi}}{k^2}
\end{eqnarray*}
for some $\varphi > 0$.
\ignore{
\begin{eqnarray*}
2+\xi-\epsilon-2\gamma'-2\theta
< 2+\xi - \frac{8+2\xi}{3} +2\delta -\epsilon \\
=\frac{\xi - 2}{3} + 2\delta -\epsilon < 0.
\end{eqnarray*}
}
Thus
\[
P(F)\leq CM^{-\varphi}\sum_{k=1}^\infty k^{-2}
\]
for 
\begin{eqnarray*}
F=
\left\{\exists_{k\in\{1,2,\frac{M^{1-\theta}}{2}\}}
\exists_{u\in[u_{k+1},u_{k+2}]} I(u)>I(u_{k})+\kappa |D(u_{k})|/2
\right\}.
\end{eqnarray*}

Similarly, for each $k\geq 2$,
\begin{eqnarray*}
P\left(
\exists_{u\in[u_{k-1},u_{k}]} \left|I(u)-I(u_{k-1})\right| \geq \frac 14\kappa |D(u_{k})|
\right)
\leq \frac{4\kappa M^\alpha }{\kappa^2 D(u_{k})^2}
\leq\frac{4M^{-\varphi}}{(k-1)^2},
\end{eqnarray*}

and thus 
\[
P(B) \leq 4M^{-\varphi}\sum_{k=1}^\infty k^{-2}
\]
for 
\[
B=
\left\{\exists_{k\in\{2,3,\frac{M^{1-\theta}}{2}\}}
\exists_{u\in[u_{k-1},u_{k}]} \left|I(u)-I(u_{k-1})\right| \geq \frac 14\kappa |D(u_{k})|
\right\}.
\]

On the event $B^c\cap F^c$, for every $k$ and every $u\in[u_{k-1},u_k]$ and $u'\in[u_{k+1},u_{k+2}]$,
we have that 
\[
z(u') < z(u_k) - k\kappa M^{\theta - 1}/2 
\ \ \ \mbox {and} \ \ \ 
z(u) > z(u_k) -k\kappa M^{\theta - 1}/4.
\]
Therefore, $z(u')<z(u)-\frac k4\kappa M^{\theta - 1} < z(u)-kM^\alpha$, and therefore, since the particle mass is bounded by $M^\alpha$, no excursion of length greater than or equal to $M^{\gamma'}$ can start at any point $u$ between $u_0$ and $u_0+\frac M2$.
By standard estimates for the size biased sequence, the probability that there is an excursion of length $M^\gamma$ starting after $u_0+M/2$ and no excursion of length larger than $M^{\gamma'}$ between $u_0$ and $u_0$ and $u_0+M/2$ decays like $\exp(-M^{\gamma-\gamma'})$. Therefore, with probability at least $1-CM^{\varphi}$, there is no such excursion starting after $u_0$.

\ignore{
\[
\gamma'+\theta - 1 > \alpha?
\]
\begin{eqnarray*}
\gamma'+\theta - 1 - \alpha > 2\gamma' - 1 - \delta -\alpha \\
= 2\gamma' - 1 - \delta - \xi + \gamma' + \epsilon
= 3\gamma' - 1 - \xi +\epsilon - \delta\\
> 3\gamma' - 1 - \xi  > \frac{4+\xi}{2} - 1 - \xi = 1-\frac \xi 2 > 0.
\end{eqnarray*}
}


\end{proof}

\begin{proof}[Proof of Lemma \ref{lem:aldous}]
Lemma \ref{lem:aldous} now follows from Lemma \ref{lem:first_block} and Lemma \ref{lem:other_blocks}.
\end{proof}

\section*{Acknowledgment}

{}First, I thank Yuval Peres and Itai Benjamini for presenting these problems
to me and for helping me during the research. I also wish to thank Omer Angel
and Elchanan Mossel for helpful suggestions. I thank Jeff Steif for his help
in improving the exposition of the paper and for presenting
to me the question leading to Theorem \ref{jeff}. I thank Michael Aizenman for
useful and interesting discussions.

I thank Mario W\"utrich for finding a mistake in an earlier version.

\noindent
Noam Berger, \\
Department of Statistics, \\
367 Evans Hall \#3860, \\
University of California Berkeley, \\
CA 94720-3860 \\
e-mail:noam@stat.berkeley.edu\\

\end{document}